\documentclass[12pt, color]{amsart}
\newcommand{\pof}{\noindent{\em Proof: }}
\usepackage{color}
\usepackage[normalem]{ulem}

\def\qdet{{\operatorname{det}}_q}

\newcommand{\s}[1]{\mathcal{#1}}

\newtheorem{Thm}{Theorem}[section]
\newtheorem{Def}[Thm]{Definition} \newtheorem{Rem}[Thm]{Remark}
\newtheorem{Lem}[Thm]{Lemma} \newtheorem{Cor}[Thm]{Corollary}
\newtheorem{Prop}[Thm]{Proposition}

\numberwithin{equation}{section}

\fontsize{.5cm}{.5cm}\selectfont\sf

\medskip

\begin{document}

\title
[Double-partitions]{Double-partition Quantum Cluster Algebras}

\author{Hans Plesner jakobsen, Hechun Zhang }
\address{
  Department of Mathematical Sciences\\ University of
Copenhagen\\Universitetsparken 5\\
   DK-2100, Copenhagen,
  Denmark\\Department of Mathematical Sciences\\Tsinghua University,
  Beijing, 100084, P. R. China} \email{jakobsen@math.ku.dk,
hzhang@math.tsinghua.edu.cn}\date{\today}

\begin{abstract} A family of quantum cluster algebras is introduced and
  studied. In general, these algebras are new, but sub-classes have been studied
previously by other authors. The algebras are indexed by double partitions or
double
  flag varieties. Equivalently, they are indexed by broken lines $L$. By grouping
  together neighboring mutations into quantum line mutations we can mutate  from
the cluster algebra of one broken line to another. Compatible pairs can be
written down. The algebras are equal to their upper cluster algebras. The
variables of the quantum seeds are given by elements of the dual canonical
basis.\end{abstract}
\subjclass[2000]{MSC 17B37 (primary),\ MSC 20G42 (primary),\ MSC 13F60
(secondary), \and MSC 32M15 (secondary)}
\keywords{Quantum algebra \and Quantized matrix algebras}
\thanks{The second author is partially supported by NSF of 
China}
 \maketitle

\section{Introduction}

A cluster algebra, as invented by Fomin and  Zelevinsky, is a
commutative algebra generated by a family of generators called
cluster variables. The generators are grouped into clusters and the
cluster variables can be computed recursively from the initial
cluster.

The theory of cluster algebras is related to a wide range of
subjects such as Poisson geometry,
 integrable systems,
 higher Teichm\"uller spaces,
 combinatorics,
 commutative and non-commutative algebraic geometry,
and the representation theory of quivers and finite-dimensional
algebras.

In \cite{Scott},  it is proved that the coordinate rings of
$SL(n,{\mathbb C})$ and its maximal double Bruhat cell
$SL(n,{\mathbb C})^{w_0,w_0}$ are cluster algebras. This is generalized in 
the recent work  \cite{CA3} where it is proved that the coordinate ring
of any double Bruhat cell $G^{u,v}$ of any semi-simple algebraic
group is a cluster algebra.

Quantum cluster algebras were introduced and studied by Berenstein and
Zelevinsky \cite{bz}. A main motivation was to understand the
dual canonical basis. Following Lusztig (\cite{lu1}), the dual canonical basis
for the coordinate
algebra ${\mathcal O}_q(M(n))$ of an $n\times n$ quantum matrix was
shown to exist in \cite{jz3}. This construction can  be carried over
to the algebra ${\mathcal O}_q(M(m,n))$ for all $m$ and $n$
verbatim.

 A quantum mutation is governed by a pair of
matrices, called a compatible pair, with certain favorable properties.
To construct a quantum cluster, one of the main difficulties is to
construct the compatible pairs. In the present paper we construct a
family of quasi commuting quantum minors of the algebra ${\mathcal
O}_q(M(m,n))$ associated to each so-called broken line $L$, and construct a
corresponding compatible pair $(\Lambda_L, B_L)$. 

The set of broken lines has a natural partial ordering with unique biggest
and smallest elements.

Let us be more specific: A {\bf 
broken line} from $(1,n)$ to $(m,1)$  is a path in ${\mathbb N}\times {\mathbb
N}$ starting at $(1,n)$  and terminating at $(m,1)$ while alternating between
horizontal and vertical segments and passing through smaller column numbers (in
the horizontal
direction)  and bigger row numbers (in the vertical direction). To each broken line
we construct in Section~\ref{6} a family of $nm$ $q$-commuting quantum minors. With the line
fixed, each of these quantum minors is uniquely given by a point $(i,j)\in
{\mathbb N}\times {\mathbb N}$ with $1\leq i\leq m$ and $1\leq j\leq n$. The
quantum cluster algebra ${\mathcal A}^-_L$ is then determined by the quantum
minors corresponding to the points on or below the line $L$. We prove that
monomials in these are members of the dual canonical basis. We introduce the
natural ordering on the set of broken lines and introduce some natural
sub-algebras. One such is ${\mathcal
    O}_q(T_L\cup L)$ which denotes the sub-algebra of ${\mathcal O}_q(M(m,n))$
  generated by the standard elements $Z_{i,j}$ of ${\mathcal O}_q(M(m,n))$
(cf. Section~\ref{333}) for which $(i,j)$ is on the line, or below it.

We introduce a special class of mutations that are called quantum line
mutations. To each triple of broken lines  $L_a,L_b$ with  $L_a,L_b\leq L$ we
can mutate by quantum line mutations from $L_a$ to $L_b$.

For a compatible pair $(\Lambda_L,B_L)$
connected with ${\mathcal
  A}^-_L$ we show that we can mutate by quantum line mutations to a bigger line
$L_1$ inside ${\mathcal
O}_q(M(m,n))$ and, by carefully keeping track, construct a  compatible
pair $(\Lambda_{L_1},B_{L_1})$ connected with  ${\mathcal
  A}^-_{L_1}$ in the process. Starting at a particularly simple broken
line, namely the one corresponding to the smallest  broken line $L^-$,  we
can, by repeated quantum line mutations, construct a compatible pair for
${\mathcal
  A}^-_L$. Thus, we   obtain
compatible pairs for all broken lines. At first they are just compatible pairs
for the smaller algebras. The algebra   ${\mathcal
O}_q(M(m,n))$ corresponds to the unique maximal broken line $L^+$. However,
mutating in the opposite direction, we get a compatible pair for this  bigger
algebra for any line. Or, indeed, mutating backwards from any bigger line
algebra to a smaller, we get a quantum seed ${\mathcal Q}_{L_1,L}$ 
for the bigger line algebra $L$
indexed by the smaller line algebra $L_1$.

Instances of such algebras have been studied in \cite{lr}, \cite{llr}, and
\cite{ga}.

The main technical result is the following: Let $A$ be an $n\times n$
matrix whose entries are non-negative integers and let $b(A)$ be the element of 
the
dual canonical basis of ${\mathcal O}_q(M(m,n))$ corresponding to this. Let
$\qdet$ denote
the quantum determinant. If $I$ denotes the $n\times n$ identity
matrix then $$\qdet=b(I)$$ and (this is  (\ref{deteq}))
$$b(A)\qdet=b(A+I).$$
Once this has been established, it can be generalized to several other
configurations involving quantum minors.

After this introduction, the article continues in Section~2 with a review of
quantum cluster algebras, followed in Section~~3 by basic facts and structures
relating to the quantized matrix algebra.  The matters concerning the technical
result (\ref{deteq}) and its generalizations, take up Sections~4 and 5.

In Section~6, using the results on dual canonical bases, we strengthen a result
of Parshall and Wang considerably. In so doing, we obtain a crucial commutation identity;
Theorem~\ref{4.4}. This result then makes it possible to introduce the class of
mutations called quantum line mutations. We also include an
observation relating this to totally positive matrices.

  In Section~7, we construct compatible pairs $(\Lambda_L^0,B^0_L)$ and 
$(\Lambda_L,B_L)$. At first just
for the algebra ${\mathcal O}_q(T_L\cup L)$, but later also for the full
algebra ${\mathcal O}_q((M(m,n))$. 
    
    Finally, in Section~8, we extend slightly a result of Goodearl and Lenagan
(\cite{gl}) saying that the $q$-determinantal ideal is prime. We then use
quantum line mutations to give an inductive proof of the following,
where
${\mathcal C}_L^-$ are the non-mutable (covariant) elements, and ${\mathcal
U}^-_L$ is the upper cluster algebra:

\smallskip

\noindent{\bf Theorem} {\em Let ${\mathcal C}_L^-=\{Y_1,\dots,Y_s\}$. Then,
$${\mathcal U}^-_L={\mathcal O}_q(T_L\cup
L)[Y_1^{\pm1},\dots,Y_s^{\pm1}]={\mathcal A}^-_L.$$}

This result is Theorem~\ref{mainthm}. As a consequence, we conclude
that in the case of ${\mathcal O}_q(M(m,n))$, the quantum  cluster algebra is
equal to its upper cluster algebra.

\subsubsection*{Acknowledgement} We would like to express our gratitude  to the referee for many helpful and clarifying suggestions.

\bigskip

\section{Basics of Quantum Cluster Algebras}

\label{2}

Throughout the paper, the base field is $K={\mathbb Q}(q)$, where
$q$ is an indeterminate over the rational numbers. To avoid terms involving
$q^{\frac12}$, we work with the square root of the $q$ used by Berenstein and
Zelevinsky; $q^2=q^2_{our}=q_{BZ}$.

 Given an integral
skew-symmetric matrix $\Lambda=(\lambda_{ij})\in M_{\mathfrak m}({\mathbb Z})$,
the  quasi polynomial algebra ${\mathcal L}(\Lambda)$
associated to the matrix $\Lambda$ is an associative algebra
generated by $x_1,x_2,\cdots, x_{\mathfrak m};
x_1^{-1},x_2^{-1},\cdots,x_{\mathfrak m}^{-1}$ with the defining relations
\begin{equation}\label{com3}x_ix_j=q^{2\lambda_{ij}}x_jx_i.
\end{equation}
Conversely, given such relations, the matrix $\Lambda=(\lambda_{ij})\in
M_{\mathfrak m}({\mathbb Z})$ will be called the $\Lambda$-matrix
of the variables $x_1\cdots,x_{\mathfrak m}$.

The set of ordered monomials
\[\{x^{\underline{a}}:=x_1^{a_1}x_2^{a_2}\cdots
x_{\mathfrak m}^{a_{\mathfrak m}}\mid\underline{a}=(a_1,a_2,\cdots,a_{\mathfrak
m})\in {\mathbb Z}^{\mathfrak m}\}\]
is a basis of ${\mathcal L}(\Lambda)$. It is well known that
${\mathcal L}(\Lambda)$ is a Noetherian domain and one can talk
about its  skew field of fractions which is  denoted  by ${\mathcal
F}(\Lambda)$. Using $\Lambda$, one can define a bilinear form on
${\mathbb Z}^{\mathfrak m}$ as follows:

\begin{eqnarray}\Lambda: {\mathbb
Z}^{\mathfrak m}\times {\mathbb Z}^{\mathfrak m}\longrightarrow {\mathbb
Z}\\\nonumber
\Lambda(\underline{a}, \underline{b})=\underline{a}\Lambda
\underline{b}^T.\end{eqnarray}

For any $\underline{a}\in {\mathbb Z}^{\mathfrak m}$, the normalized monomial is
defined as
\[x(\underline{a})=q^{\sum_{i<j}\lambda_{ji}a_ia_j}x^{\underline{a}}.\]
The map \begin{equation}\label{q-bar} \forall i=1,\dots, {\mathfrak m}:
x_i\mapsto x_i,\ q\mapsto q^{-1}
\end{equation}
extends to a ${\mathbb Q}$-algebra anti-automorphism, denoted by $\ell\mapsto\overline\ell$, which actually does not
depend on the ordering. Then
\begin{equation}\overline{x(\underline{a})}=x(\underline{a}).\label{bar-2}
\end{equation}

It is easy to check that
\[x(\underline{a})x(\underline{b})=q^{\Lambda(\underline{a},
\underline{b})}x(\underline{a}+\underline{b}),\]
which, of course, is equivalent to the commutation relations (\ref{com3}).

Denote by $K^*:={\mathbb Q}(q)-\{0\}$ the multiplicative group of
non-zero elements. The group $(K^*)^{\mathfrak m}$ acts on ${\mathcal
L}(\Lambda)$ as an automorphism group. Explicitly, for any
$\underline{h}=(h_1,h_2,\cdots, h_{\mathfrak m})\in (K^*)^{\mathfrak m}$, it
acts
on
${\mathcal L}(\Lambda)$ according to the formulae
\[\underline{h}(x_i)=h_ix_i\text{ for all }i.\]

\begin{Rem} If a subspace $S\subset {\mathcal
A}(\Lambda)$ is invariant under the action of the group $(K^*)^{\mathfrak m}$,
then it is spanned by the monomials that it contains.\end{Rem}

\medskip

\medskip

In \cite{bz}, the notion of a quantum cluster algebra was introduced. Let us
recall the definition.

\begin{Def}
\label{def:compatible-triple} Let $B$ be an ${\mathfrak m} \times {\mathfrak n}$
 integer
matrix with rows labeled by $[1,\mathfrak m]$ and columns labeled by an
$\mathfrak n$-element subset $ex \subset [1,\mathfrak m]$.
Let $\Lambda$ 
be a skew-symmetric
$\mathfrak m \times \mathfrak m$ integer matrix with rows and columns labeled
by $[1,\mathfrak m]$. We say that a pair $(\Lambda,  B)$ is \emph{compatible}
if,  for every $j \in ex$ and $i \in [1,\mathfrak m]$, we have
$$\sum_{k = 1}^{\mathfrak m} b_{kj} \lambda_{ki} = \delta_{ij} d_j $$
for some positive integers $d_j \,\, (j \in ex)$. The $\mathfrak n \times
\mathfrak n$ sub-matrix  of $B$ corresponding to the subset $ex$ is called the
principal part of $B$. {{} We insist  throughout this article, that $\forall j: d_j=2$.}
\end{Def}

\noindent If one arranges the symbols such that $ex=\{1,2,\dots,\mathfrak
n\}$,  the compatibility condition states that 
the $\mathfrak n \times \mathfrak m$ matrix $\tilde D = B^T \Lambda$ consists of
the two
blocks: the $\mathfrak n \times \mathfrak n$ diagonal matrix $D$ with positive
integer
diagonal entries $d_j$, and the $\mathfrak n \times (\mathfrak m-\mathfrak n)$
zero block.

With the above setup, the triple $(\{x_1, x_2,\cdots, x_{\mathfrak m}\},
\Lambda, B)$ is an example of a {\bf quantum seed} of ${\mathcal F}(\Lambda)$.
The notion of a quantum seed is more general than the one presented here, but
ours suffices for the purposes below.
The variables $x_i$ are called {\em quantum cluster variables}. The variables
$x_i, i\in ex$ are called {\em mutable variables} and the set of these is called {\em the cluster}. The variables
$x_j, j\notin ex$ are called {\em non-mutable variables}.

\medskip

{Notice that if $\underline{a}=(a_1,a_2,\dots,a_{\mathfrak m})$ and
$\underline{f}=(f_1,f_2,\dots,f_{\mathfrak m})$ are vectors then
\begin{Lem}\label{2.22}\begin{equation}
\Lambda(\underline{a})^T=(\underline{f}
)^T\Leftrightarrow\forall
i:x_ix^{\underline{a}}=q^{2f_i}x^{\underline{a}}x_i.\end{equation}
In particular, if there exists a $j$ such that $\forall i: f_i=-\delta_{i,j}$
then the column vector $\underline{a}$ can be the $j$th column in the matrix $B$
of a compatible pair.
\end{Lem}
\noindent However simple this actually is, it will have a great importance later
on.}

\medskip

Denote by  $e_1,e_2,\cdots,e_{\mathfrak m}$ the standard basis of ${\mathbb
Z}^{\mathfrak m}$. For a given compatible pair $(\Lambda, B=(b_{ki}))$, one can
{\em mutate} the cluster in the direction of  $i\in ex$, thereby
obtaining a new cluster
whose variables are $x_1,\cdots, x_{i-1}, x_i^\prime,x_{i+1}$, $\cdots$, $
x_{\mathfrak m}$. The unique new variable is defined by
\begin{equation}x_i^\prime
=x(\sum_{b_{ki}>0}b_{ki}e_k-e_i)+x(\sum_{b_{ki}<0}-b_{ki}e_k-e_i).\label{muta}\end{equation}
One can check that $x_1,\cdots, x_{i-1}, x_i^\prime,x_{i+1},\cdots,
x_{\mathfrak m} $ is a q-commuting family.

We will extend matrix mutations to those of compatible pairs. Fix an
index $i \in ex$. The matrix $B_i' = \mu_i(B)$ can be
written as
\begin{equation}
\label{eq:mutation-product}B_i' = E_i\,B\, F_i\, ,
\end{equation}
where

\noindent $\bullet$ $E_i$
is the $\mathfrak m\times \mathfrak m$ matrix with entries
\begin{equation}
\label{eq:E-entries} e_{ab} =\left\{
\begin{array}{lc}
\delta_{ab} & \textrm{if } b \neq i;\\[.05in]
- 1 & \textrm{if }a = b = i; \\[.05in]
\max(0, -b_{ai})& \textrm{if } a \neq b = i.
\end{array}\right. 
\end{equation}

\noindent $\bullet$ $F_i$
is the $\mathfrak n\times \mathfrak n$ matrix with rows and columns labeled by
$ex$, and
entries given by
\begin{equation}
\label{Zel} f_{ab} =\left\{
\begin{array}{lc}
\delta_{ab} & \textrm{if } a \neq i; \\[.05in]
- 1 & \textrm{if } a = b = i; \\[.05in]
\max(0,  b_{ib})& \textrm{if } a = i \neq b. 
\end{array}\right.
\end{equation}

The triple $(\{x_1,\cdots, x_{i-1}, x_i^\prime, x_{i+1},\cdots,
x_{\mathfrak m}\}, \Lambda_i=E_i^T\Lambda E_i, B'_i)$ is also a quantum
seed. The above process of passing  from a quantum seed to another is called a
quantum mutation in the direction $i$. We say that two quantum seeds
are mutation equivalent if they can be obtained from each other by a
sequence of quantum mutations. In the general definition of \cite{bz} there is an
additional parameter $\varepsilon=\pm1$ in the definition of the matrices
$E_i,F_i$. Throughout this article we restrict to $\varepsilon=1$ and for this
reason we suppress it. 

Given a quantum seed, let ${\mathcal S}$ be the
set of all quantum
seeds which are mutation equivalent to the given one.  The
quantum cluster algebra ${\mathcal A}({\mathcal S})$ associated to the given
quantum seed is the
${\mathbb Q}(q)$ sub-algebra of ${\mathcal F}(\Lambda)$  generated by
all quantum cluster variables contained in ${\mathcal S}$.

\section{The quantum matrices and the dual canonical basis}

\label{333}

The coordinate algebra  ${\mathcal O}_q(M(m,n))$ of the quantum
$m\times n$ matrix is an associative algebra, generated by elements
$Z_{ij},i=1,2,\cdots, m;j=1,2,\cdots,n$, subject to the following
defining relations:
\begin{eqnarray}\label{relations1}Z_{ij}Z_{ik}&=&q^2Z_{ik}Z_{ij} \text{ if }
j<k,\\
 Z_{ij}Z_{kj}&=&q^2Z_{kj}Z_{ij} \text{ if }i<k,\\
Z_{ij}Z_{st}&=&Z_{st}Z_{ij}\text{ if } i>s, j<t,\\
Z_{ij}Z_{st}&=&Z_{st}Z_{ij}+(q^2-q^{-2})Z_{it}Z_{sj} \text{ if
}\label{relations4}\label{344}
i<s, j<t.\end{eqnarray} 

The associated quasi-polynomial algebra ${\mathit o}_q(M(m,n))$ of the quantum
$m\times n$ matrix is an associative algebra, generated by elements
$z_{ij},i=1,2,\cdots, m;j=1,2,\cdots,n$, subject to the following
defining relations:
\begin{eqnarray}\label{relations2}z_{ij}z_{ik}&=&q^2z_{ik}z_{ij} \text{ if }
j<k,\\
 z_{ij}z_{kj}&=&q^2z_{kj}z_{ij} \text{ if }i<k,\\
z_{ij}z_{st}&=&z_{st}z_{ij} \text{ in all other
cases}.\label{3.7}\end{eqnarray}

For any matrix $A=(a_{ij})_{1\le i\le m,1\le
j\le n}\in M_{m,n}({\mathbb Z}_+)$, where  ${\mathbb Z}_+=\{0,1,\cdots\}$,
we define a monomial $Z^A$  by
\begin{equation} Z^A=\Pi_{i,j=1}^nZ_{ij}^{a_{ij}},\end{equation}
where the factors are arranged in the descending lexicographic order on
$I(m,n)=\{(i,j)\mid i=1,2,\cdots,m; j=1,\cdots,n\}$ given by
$(1,1)>(1,2)>\dots>(1,n)>(2,1)>\dots$. We define  similar elements
$z^A\in{\mathit o}_q$. It is well
known that the set $\{Z^A|A\in M_{m,n}({\mathbb Z}_+)\}$ is a basis
of the algebra ${\mathcal O}_q(M(m,n))$.

From the defining relations (\ref{relations1}) - (\ref{3.7}) of the
algebras
${\mathcal O}_q(M(m,n))$ and ${\mathit o}_q(M(m,n))$ it is easy to show the
following lemma. The last statement in the lemma, though trivial,  is included
for its usefulness.

\begin{Lem} \label{3.1}
The mapping
\begin{eqnarray}^-:Z_{ij}&\mapsto &Z_{ij}\\\nonumber
q&\mapsto &q^{-1}\end{eqnarray} extends to an 
anti-automorphism of the algebra ${\mathcal O}_q(M(m,n))$ as an
algebra over ${\mathbb Q}$. The mapping \begin{eqnarray}^-:z_{ij}&\mapsto
&z_{ij}\\\nonumber
q&\mapsto &q^{-1}\end{eqnarray} extends to an 
anti-automorphism of the algebra ${\mathit o}_q(M(m,n))$ as an
algebra over ${\mathbb Q}$.

There is an obvious anti-automorphism of  the tensor algebra over the vector space $M(m,n)$ of which the given anti-automorphism  of ${\mathcal O}_q(M(m,n))$ 
is the quotient map.

 A similar statement holds in   ${\mathit o}_q(M(m,n))$.
\end{Lem}

For any $A=(a_{ij})\in M_{m,n}({\mathbb Z}_+)$,
$$\underline{ro}(A):=(\sum_j
a_{1j},\cdots,\sum_ja_{mj}):=(\underline{ro}_1,\underline{ro}_2,\cdots,
\underline{ro}_m).$$
This is called the row sum of $A$. 
$$\underline{co}(A):=(\sum_j a_{j1},\cdots,\sum_j
a_{jn}):=(\underline{co}_1,\underline{co}_2,\cdots,\underline{co}_n).$$This is
called the column sum of $A$.

\medskip

The following result follows easily from the defining relations 
(\ref{relations1}) - (\ref{relations4}):
\begin{Lem}Let
\begin{displaymath}
  Z^AZ^B=\sum_Ca^{A,B}_CZ^C 
\end{displaymath}
where $a^{A,B}_C\in{\mathbb Z}[q^2,q^{-2}]$. Then $\forall a^{A,B}_C\neq0$:
\begin{eqnarray*}
  \underline{ro}(C)&=&\underline{ro}(A)+\underline{ro}(B)\\
\underline{co}(C)&=&\underline{co}(A)+\underline{co}(B).
\end{eqnarray*}
\end{Lem}

\medskip

From the defining relations we also have

\begin{equation}\label{bar}\overline{Z^A}=E(A)Z^A+\sum_{B<A}c_B(A)Z^B,
\end{equation}
where
$$E(A)=q^{-2(\sum_i\sum_{j>k}a_{ij}a_{ik}+\sum_i\sum_{j>k}a_{ji}a_{ki})}$$
and $\forall B<A$: $c_B(A)\neq0\Rightarrow \underline{ro}(B)=\underline{ro}(A)$,
and  $\underline{co}(B)=\underline{co}(A)$. Here, $c_B(A)\in{\mathbb
Z}[q^2,q^{-2}]$, and the lexicographic order on $M_{m,n}({\mathbb Z}_+)$, obtained
by augmenting the previous order on $I(m,n)$ by the natural order on ${\mathbb
Z}_+$,  is denoted $\leq$.

\medskip

Let \begin{equation}
	N(A)=q^{-\sum_i \sum_{j>k}a_{ij}a_{ik}-\sum_i
\sum_{j>k}a_{ji}a_{ki}}\ \textrm{ and }Z(A)=N(A)Z^A.
	\label{norma}
\end{equation}

From  (\ref{bar}) we trivially have \begin{equation}
\overline{Z(A)}=Z(A) \ \textrm{ modulo lower order terms.}
	\label{bar-3}
\end{equation} 

\noindent In lack of better words we introduce:
\begin{Def} We call $Z(A)$ the normalized form of $Z^A$. We call $N(A)$ the
normalization factor. 
\end{Def}

\medskip

Let $i<s$ and $j<t$. Set $E_{i,j,s,t}=E_{i,j}+E_{s,t}-E_{i,t}-E_{s,j}$, where
for any of the mentioned pairs $(a,b)$, $E_{a,b}$ is the $(a,b)$th matrix unit. 
Upon
rewriting $\overline{Z^A}$ according to our lexicographic order as in (\ref{bar}), one picks up
terms $c_{A'}Z^{A'}$,  where $A'$ is
obtained  from $A$ by subtraction of elements of the form $E_{i,j,s,t}$. The
next result follows directly from (\ref{norma}).  

\begin{Lem}\label{lemnorm}If $A'=A-E_{i,j,s,t}$, then $$N(A')=N(A)q^{4-2(a_{ij}+a_{st}-a_{it}-a_{sj})}.$$
\end{Lem}

\medskip

To facilitate the following proofs, we introduce a notion of a level in
$M_{m,n}({\mathbb Z}_+)$: 

\begin{Def} Let \begin{equation}{\mathcal D}=\{E_{i,j,s,t}\mid
i<s\textrm{ and }j<t\}.\end{equation} The matrix $A\in M_{m,n}({\mathbb Z}_+)$ is of level $L(A)=0$ if
there are
no elements $D\in{\mathcal D}$ and $A_1\in M_{m,n}({\mathbb Z}_+)$ such that
$A=D+A_1$. Let ${\mathcal L}_0$ denote the set of matrices of level 0. We define
the level  $L(A)$ of any $A$ not of level zero by
\begin{displaymath}
	L(A):=\max\{r\in{\mathbb N}\mid \exists D_1,\dots,D_r\in{\mathcal D},
\exists A_0\in{\mathcal L}_0: A=D_1+\dots+D_k +A_0\}.
\end{displaymath}
 (It is easy to see that this maximum is finite).
\end{Def}

Notice that if $L(A)=0$ then one can reorder $\overline{Z(A)}$ without invoking
the relation (\ref{344}). Thus, $\overline{Z(A)}={Z(A)}$.

\medskip

\begin{Lem}\label{addlem}In the equation (\ref{bar}), if $B<A$ and $c_B(A)\neq0$
then $L(B)<L(A)$.
\end{Lem}

\proof View the right hand side of (\ref{bar}) as the result of the
reordering of $\overline{Z^A}$ according to our chosen ordering. The terms with
$B<A$ must then have their origins in the application of relation (\ref{344})
at least once since otherwise we can get to $Z^A$ using solely the
other three relations. In this case,   $\overline{Z(A)}=Z(A)$. Any
application of (\ref{344}) clearly leads, modulo terms proportional to $Z^A$,
to terms of lower level. \qed

\medskip

 Set
 $$L^*=\oplus_{A\in M_{m,n}({\mathbb Z}_+)}{\mathbb Z}[q]Z(A).$$

\medskip

\begin{Prop}\label{basis}
There is a unique ${\mathbb Z}[q]$-basis $B^*=\{b(A)|A\in M_{m,n}({\mathbb
Z}_+)\}$
of $ L^*$ in which each element $b(A)$ is determined uniquely by the following
conditions:
\begin{enumerate}
\item $\overline{b(A)}=b(A)$.

\item $b(A)=Z(A)+\sum_{B<A} h_B(A)Z(B)$ where $h_B(A)\in q^2{\mathbb Z}[q^2]$
and $\underline{ro}(B)=\underline{ro}(A), \underline{co}(B)=\underline{co}(A)$.
\end{enumerate}
 The basis $B^*$ is called the dual canonical basis of ${\mathcal
O}_q(M(m,n))$.\end{Prop}

\begin{Cor}\label{changeofb}If we number our basis vectors in the two bases
${\mathcal
B}_1=\{b(B)\mid B\in M_{m,n}({\mathbb Z}_+)\}$ and ${\mathcal B}_2=\{Z(B)\mid
B\in M_{m,n}({\mathbb Z}_+)\}$ according to the lexicographic ordering then the
change of basis matrices are lower triangular with 1's in the diagonal and
elements from $q^2{\mathbb Z}[q^2]$ in all other non-zero
positions.\label{basiscor}
\end{Cor}

\proof [Proof of Proposition~\ref{basis} and Corollary~\ref{changeofb}]  Noticing the $q^2$ factors in
Lem\-ma~\ref{lemnorm}, this can be proved in analogy with Lusztig (\cite[2.
Proposition]{lu1}), see also
\cite[Theorem~3.5]{jz3}.
However, we will sketch a proof for clarity: We proceed to prove Proposition~\ref{basis} by induction on the
level $k$, utilizing that if the proposition holds up to level $k$ then
so does Corollary~\ref{basiscor}. The case of level 0 is trivial since if
$L(A)=0$ 
then $b(A)=Z(A)$. Suppose then that the result holds up to, and including level
$k$ and let $A$ be of level $k+1$. It follows from (\ref{bar}) and
Lemma~\ref{addlem} together with
Corollary~\ref{basiscor} (up to level $k$) that \begin{equation}\label{left}
\overline{Z(A)}-Z(A)=\sum_{B<A;\ L(B)<L(A)}\frac{c_B(A)}{N(A)N(B)}Z(B) =\sum_{B<A;\
L(B)<L(A)}h_Bb(B) 
\end{equation}
with elements $h_B\in {\mathbb Z}[q^{2},q^{-2}]$. Since the left hand side of
(\ref{left}) is skew under the bar operator, each $h_B$ can be decomposed as
$h_B=h_B^++h_B^-$ with $h_B^+\in q^2{\mathbb Z}[q^2]$ and
$h_B^-=-\overline{h_B^+}$. Then \begin{equation}
                                 b(A)=Z(A)+\sum_{B<A}h_B^+b(B) 
                                \end{equation}
is the unique solution. Invoking Corollary~\ref{basiscor} (up to $k$) once again, the proof
is complete. \qed

\medskip

The following simple principle is very useful:

\smallskip

\begin{Prop}\label{simple}
If $m_1\leq m$ and $n_1\leq n$ we may view ${\mathcal O}_q(M(m_1,n_1))$ as the
sub-algebra of ${\mathcal
O}_q(M(m,n))$ generated by the elements of (some) $m_1$ rows and $n_1$ columns.
If, correspondingly, we consider $M_{m_1,n_1}({\mathbb Z}_+)\subseteq 
M_{m,n}({\mathbb Z}_+)$ then for any $A\in M_{m_1,n_1}({\mathbb Z}_+)$, upon
these identifications,  the basis vector $b(A)\in {\mathcal O}_q(M(m_1,n_1))$ is
also a basis vector in ${\mathcal O}_q(M(m,n))$. 

If, under such identifications, ${\mathcal O}_q(M(m_1,n_1))$ and ${\mathcal
O}_q(M(m_2,n_2))$ are two {\sf commuting} sub-algebras of ${\mathcal
O}_q(M(m,n))$  and if $b(A_i)\in {\mathcal O}_q(M(m_i,n_i))$, $i=1,2$, are
members of the respective dual canonical bases, then $b(A_1+A_2)=b(A_1)b(A_2)$
is in the dual canonical basis of ${\mathcal O}_q(M(m,n))$.
\end{Prop}

\proof The relations, the bar operator, and the order on the sub-algebras are
restrictions of the relations, the bar operator, and the order on the full
algebra. The result then follows by the uniqueness. \qed 
                        
\begin{Rem}The commutativity condition in Proposition~\ref{simple} is equivalent
to having all canonical generators of one sub-algebra positioned NE of the
other.
\end{Rem}

\bigskip

If $m=n$,  one may define the quantum determinant $det_q$ as
follows:
\begin{eqnarray}\label{2b} \qdet(n)=\qdet&=&\Sigma_{\sigma\in
S_n}(-q^2)^{\ell(\sigma)}Z_{1,\sigma(1)}Z_{2,\sigma(2)} \cdots
Z_{n,\sigma(n)}\\\label{3b}&=&\Sigma_{\delta\in
S_n}(-q^2)^{\ell(\delta)}Z_{\delta(1),1}Z_{\delta(2),2} \cdots
Z_{\delta(n),n}. \end{eqnarray}

\medskip

We recall some results  from \cite{nmy} regarding the Quantum
Laplace Expansion:

Suppose $I=\{i_1<1_2< \dots < i_r\}$ and $J=\{j_1<j_2< \dots <
j_r\}$ are subsets of $I=\{1,2,\dots, n\}$. Define
\begin{eqnarray}\label{2} \xi^I_J&=&\Sigma_{\sigma\in
S_r}(-q^2)^{\ell(\sigma)}Z_{i_1,j_{\sigma(1)}}Z_{i_2,j_{\sigma(2)}}
\cdots Z_{i_r,j_{\sigma(r)}}\\\label{3}&=&\Sigma_{\tau\in
S_r}(-q^2)^{\ell(\tau)} Z_{i_{\tau(1)},j_1}Z_{i_{\tau(2)},j_2} \cdots
Z_{i_{\tau(r)},j_r}.\end{eqnarray} These elements are called quantum minors.
Notice that they are only defined if $\#J=\#I$. For
two subsets $I,J\subseteq\{1,2,\dots, n\}$, the symbol ${sgn}_q(I;J)$ is
defined by
\begin{equation}
  \label{eq:oe}
  {sgn}_q(I;J)=\left\{
      \begin{array}{ll}0&\textrm{ if }I\cap J\neq\emptyset\\
(-q^2)^{\ell(I;J)}&\textrm{ if }I\cap J=\emptyset,
      \end{array}\right.
\end{equation}
where $\ell(I;J)=\#\{(i,j)\mid i\in I, j\in J,i>j\}$. Then,

\begin{eqnarray}\label{A}
  {Sgn}_q(J_1;J_2)\xi^{I}_{J}&=&\sum_{I_1\cup
    I_2=I}\xi^{I_1}_{J_1}\xi^{I_2}_{J_2}{Sgn}_q(I_1;I_2)\\
{Sgn}_q(J_1,J_2)\xi^{J}_{I}&=&\sum_{I_1\cup
I_2=I}\xi^{J_1}_{I_1}\xi^{J_2}_{I_2}{Sgn}_q(I_1;I_2)\label{B}
\end{eqnarray}

 If $m=n$ and
$I=\{1,2,\cdots,n\}\setminus\{i\},J=\{1,2,\cdots,n\}\setminus\{j\}$,
$\xi^I_J$ will  (occasionally) be denoted by $A(i,j)$.

The following was proved by Parshall and Wang in
\cite{pw}:

\begin{Prop}$\qdet$ is central. Furthermore, let $i,k\le n$ be fixed integers.
Then

\begin{equation}\label{325}\delta_{i,k}\qdet=\sum_{j=1}^n(-q^2)^{j-k}Z_{i,j}A(k,
j)=\sum_j(-q^2)^{i-j}
A(i,j)Z_{k,j}\end{equation}
\begin{equation}=\sum_j(-q^2)^{j-k}Z_{j,i}A(j,k) =
  \sum_j(-q^2)^{i-j}A(j,i)Z_{j,k}.\end{equation} 
\end{Prop}

It is of key importance for the rest of the article to note the following which
is proved by an easy induction argument using (\ref{325}) while invoking the
uniqueness of the dual canonical basis:

\begin{Cor}
$$\overline{\qdet}=\qdet = b(I) .$$
Thus, all quantum minors are members of the dual canonical basis.
\end{Cor}

\medskip

 \begin{Def} An element $x\in {\s O}_q(M(m,n))$ is called covariant if for any
$Z_{ij}$ there exists an integer $n_{i,j}$
 such that \begin{equation}xZ_{i,j}=q^{2n_{i,j}}Z_{i,j}x.\end{equation} Clearly,
$Z_{1,n}$ and $Z_{m,1}$\label{covdef}
 are covariant. 
Two elements $x,y\in {\s O}_q(M(m,n))$ are said to $q$-commute if there exists
an integer $p$ such that 
$$xy=q^{2p}yx.$$

\end{Def}   \medskip

 Let ${\det}_q(t)={\xi}^{\{1,\cdots,t\}}_{\{n-t+1,\cdots,n\}}$, for
$t=1,2,\cdots,min\{m,n\}$.
 It is easy to extend \cite[Theorem~4.3]{jz1} from the $n\times n$ case to the
general rectangular case:  \smallskip
 
 \begin{Prop}\label{cminor}The element ${\det}_q(t)$ is covariant for all $t$.
 More precisely, let $M_t^-=\{(i,j)\in {\mathbb N}^2\mid 1\le i\le t\text{ and
}1\le j\le n-t\}$, $M_t^+=\{(i,j)\in
 {\mathbb N}^2\mid t+1\le i\le m\text{ and }n-t+1\le j\le n\}$,
$M_t^l=\{(i,j)\in {\mathbb N}^2\mid t+1\le i\le
 m
 \text{ and }1\le j\le n-t \}$, and $M_t^r=\{(i,j)\in {\mathbb N}^2\mid 1\le
i\le t\text{ and }n-t+1\le j\le n\}$.
 \begin{eqnarray} Z_{i,j}{\det}_q(t)&=&{\det}_q(t)Z_{i,j}
 \text{ if }(i,j)\in M_t^l\cup M^r_t,\\\nonumber
Z_{i,j}{\det}_q(t)&=&q^2{\det}_q(t)Z_{i,j}
 \text{ if }  (i,j)\in M_t^-,\text{ and }\\\nonumber
Z_{i,j}{\det}_q(t)&=&q^{-2}{\det}_q(t)Z_{i,j}
 \text{ if }(i,j)\in M^+_t. \end{eqnarray}
\end{Prop}

\medskip

Recall from \cite[Lemma~3.3]{jz1} the result for quantum $2\times2$ matrices:

\begin{equation}
 \forall a\in{\mathbb N}:\
Z_{2,2}^aZ_{1,1}=Z_{1,1}Z_{2,2}^a+q^{-2}(1-q^{4a})Z_{2,2}^{a-1}Z_{2,1}Z_{1,2}.
\end{equation}

\medskip

For later purposes, we need the following results for $n\times n$ matrices
regarding $Z_{n,n}A(n,n)$: Using (\ref{325}) we can define elements $M_1,M_2$ by
\begin{eqnarray*}
  \qdet&=&\sum_{j=1}^n(-q^2)^{j-n}Z_{n,j}A(n,j)=\sum_{j=1}^{n}(-q^2)^{n-j}
A(n,j)Z_{n,j}\\&=&Z_{n,n}A(n,n)+M_1\\
&=&A(n,n)Z_{n,n}+M_2.
\end{eqnarray*}

If one removes the variables in the $j$th row while adding a new $0$th row, an easy application of Proposition~\ref{cminor} gives that for $j\neq n$,
$Z_{n,n}A(n,j)=q^{-2}A(n,j)Z_{n,n}$ and similarly for $A(j,n)$, and it then
follows that
$$Z_{n,n}M_i=q^{-4}M_iZ_{n,n}\textrm{ for }i=1,2.$$ (This result also follows
from \cite[Lemma~4.5.1]{pw}.) In the ring of fractions of ${\mathcal O}_q(M(n))$
we can write
$A(n,n)=(Z_{n,n})^{-1}\left(\qdet
  -M_1\right)=\left(\qdet
  -M_2\right)(Z_{n,n})^{-1} $, which is useful since, by what we have just proved, the terms $M_1,M_2$  have simple $q$-relations with
$Z_{n,n}$. Thus, 

\begin{eqnarray*}
  Z_{n,n}A(n,n)&=&\qdet-M_1=\qdet -q^{-4}M_2\\
A(n,n)Z_{n,n}&=&\qdet-q^4M_1=\qdet-M_2, \textrm{ and hence }\\
\left[Z_{n,n},A(n,n)\right]&=&(q^4-1)M_1\\&=&(1-q^{-4})M_2.
\end{eqnarray*}

Notice that all monomials in $M_2$ contain factors of $q^{2\ell}$ with
$\ell\geq1$.

\smallskip More generally, it follows by induction that for all $r\in{\mathbb N}$,
\begin{eqnarray*}\label{commu0}
  \left[Z_{n,n}^r,A(n,n)\right]&=&q^4(1-q^{-4r})M_1Z_{n,n}^{r-1}\\
&=&(1-q^{-4r})M_2Z_{n,n}^{r-1}.
\end{eqnarray*}

\noindent Likewise, for all $r\in{\mathbb N}$,

\begin{eqnarray}
  \left[Z_{1,1}^r,A(1,1)\right]&=&-(1-q^{-4r})Z_{1,1}^{r-1}N_1,\textrm{ where}\\
N_1&=&\sum_{j=2}^n(-q^2)^{j-1}Z_{1,j}A(1,j)\label{commu}\\
&=&\sum_{\sigma\in S_n; \sigma(1)\neq
1}(-q^2)^{\ell(\sigma)}Z_{1,\sigma(1)}Z_{2,\sigma(2)} \cdots
Z_{n,\sigma(n)},\label{commu2}
\end{eqnarray}
where, for each $\sigma$ in the last sum, $\ell(\sigma)\in{\mathbb N}$. This
observation is an important ingredient in the proof of Theorem~\ref{thm} below.
\label{somecom}
\medskip

\section{${\det}_q$ and dual canonical bases}

In this section, $n=m$ throughout. The following result is of key importance.
However simple to formulate, it is
remarkably difficult to prove.

\bigskip
\begin{Thm}\label{thm} For all $A\in M_n({\mathbb Z}_+)$ there are integers
$c_B\in\{-1,0,1\}$ and integers $\gamma_B>0$ such that
  \begin{equation}\label{det-for}
  Z(A)\cdot{\det}_q=Z(A+I)+ \sum_{B<A+I} q^{2\gamma_B}c_BZ(B),
\end{equation}
where each $B$ furthermore has the same row and column sums as $A+I$.
In particular,
\begin{equation}\label{deteq}
  b(A)\cdot{\det}_q=b(A+I)=b(A)b(I).
\end{equation}
\end{Thm}

\bigskip

\pof By using Corollary~\ref{changeofb}, the second claim follows easily from
the first.  We proceed to prove
(\ref{det-for}) by induction on the number $c$ such
that there are  non-zero
elements in at most the columns $1,\cdots,c$ of $A$. For a fixed such $c$ we
proceed by
induction on the number $r$ such that there are non-zero elements in at most the
rows
$1,\dots,r$ in the $c$th column. Notice that the 
formula (\ref{det-for}) holds for any $A$ with non-zero entries at most
in the first column. Indeed, as follows by an elementary computation, 
\begin{equation}Z(A)Z(E_\sigma)(-q^2)^{\ell(\sigma)}=
(-q^2)^{\ell(\sigma)}q^{2(\underline{co}_1-(a_{i,1}+a_{i+1}+\dots+a_{n,1}))}
Z(A+E_\sigma),\label{column1}\end{equation}
where $E_\sigma$ is the matrix of the permutation $\sigma$, $i=\sigma^{-1}(1)$,
and $Z(E_\sigma)=Z^{E_\sigma}$.

It is likewise easy to see
that if the theorem holds for any $A$ with non-zero entries in at
most the first $c$ columns $1,2,\cdots, c$, then it is also true if
we replace $A$ by   $A+a_{1,c+1}E_{1,c+1}$ for any ${a_{1,c+1}}\in{\mathbb N}$. Here it suffices to observe that $Z(A+a_{1,c+1}E_{1,c+1})=q^{a_{1,c+1}\underline{ro}_1(A)}Z^{a_{1,c+1}E_{1,c+1}}Z(A)$. When we multiply by $\qdet$ from the left, we obtain elements of the form $q^{a_{1,c+1}\underline{ro}_1(A)}Z^{a_{1,c+1}E_{1,c+1}}Z(A')$, where $\underline{ro}(A')=(1,1,\dots,1)+\underline{ro}(A)$, and similarly for the column sums. Due to the special form of $Z^{a_{1,c+1}E_{1,c+1}}Z(A)$ it is a matter of simple bookkeeping to verify the claim here.

Now let us assume that the theorem holds up to the $r$th row in the
$c$th column, where $r<n$. Let  $Z^{A_0}$ correspond to a matrix $A_0$ fulfilling
the requirements up to, and including, row $r$ and column $c$,  and
consider $A=A_0+a_{r+1,c}\cdot E_{r+1,c}$. 

Before getting further into the details, let us remark that the two
lexicographic orderings $(1,1)>(1,2)>\dots>(1,n)>(2,1)>\dots$ and
$(1,1)>(2,1)>\dots>(n,1)>(2,1)>\dots$ have the same monomials. By this we mean
that if $Z^A$ is written according to one of the orderings, then rewriting it
according to the other will not create auxiliary terms. Indeed, not even a 
different coefficient. Let us denote the former ordering by row-column and the
latter by column-row.

Consider 
\begin{displaymath}
Z(A)\cdot   \Sigma_{\sigma\in S_n}(-q^2)^{\ell(\sigma)}
Z_{1,\sigma(1)}Z_{2,\sigma(2)} \cdots Z_{n,\sigma(n)}.
\end{displaymath}
Set
$\alpha=(\sum_{k=1}^{c-1}a_{(r+1),c}a_{(r+1),k}+\sum_{t=1}^r a_{r+1,c}a_{t,c})$.
Then,
\begin{equation}\label{reorder}
  Z(A)=Z(A_0)\cdot q^{-\alpha}
    Z_{r+1,c}^{a_{r+1,c}}.
\end{equation}

The task now is to order each summand in $$Z(A)\cdot
\Sigma_{\sigma\in S_n}(-q)^{2\ell(\sigma)}
Z_{1,\sigma(1)}Z_{2,\sigma(2)} \cdots Z_{n,\sigma(n)}$$
lexicographically. To do so, we will group the terms in $\qdet$ together
strategically into sums of products of quantum minors. {These minors will then be ordered collectively while using their $q$-commutation relations.}

\medskip

We can safely assume  $1\leq r\leq n-1$.  Decompose
$$\{(i,j)\in{\mathbb N}\mid 1\leq i,j\leq n\}=R_1\cup R_2\cup R_3 \cup R_4,$$
$R_1=\{(i,j)\mid 1\leq i\leq r;\ 1\leq j\leq c-1\}$,  $R_2=\{(i,j)\mid 1\leq
i\leq r;\ c\leq j\leq n\}$, $R_3=\{(i,j)\mid r+1\leq i\leq n;\ 1\leq j\leq
c-1\}$, and $R_4=\{(i,j)\mid r+1\leq i\leq n;\ c\leq j\leq n\}$. Consider
({\ref{A}) applied to $\qdet$  with $J_1=\{1,\dots,c-1\}$ and
$J_2=\{c,\dots,n\}$. Then apply (\ref{B}) to each $\xi^{I_i}_{J_i}$, $i=1,2$  
based on a decomposition $I_i=R^{(1)}(I_i)\cup  R^{(2)}(I_i)$ of $I_i$; $R^{(1)}(I_i)\subseteq \{1,\dots, r\}$, and 
$R^{(2)}(I_i)\subseteq \{r+1,\dots, n\}$. The result is a formula
\begin{equation}\label{det-min}
{\det}_q=\sum_{M_1,M_2,M_3,M_4}c_{M_1,M_2,M_3,M_4}M_1M_2M_3M_4,
\end{equation}
where each $c_{M_1,M_2,M_3,M_4}$ is $\pm(q^2)^p$ for some non-negative integer
$p$, and each $M_i$,
$i=1,2,3,4$, is a quantum minor with entries from $R_i$, $i=1,\dots,4$. Notice
that it
follows from the defining relations that $M_2M_3=M_3M_2$. Not all combinations
of quantum minors will occur with a non-zero coefficient, of course. For
instance, no pair
can share a row or a column. 

Let ${\mathcal R}_{1,2,3}$ denote the set of matrices over ${\mathbb Z}_+$ with
non-zero entries at
most in the positions of $R_1\cup R_2\cup R_3$,  let ${\mathcal R}_{4}$ denote
the set of matrices with non-zero entries at most in the positions of $R_4$, and
let ${\mathcal M}_4$ denote the set of quantum minors having entries from $R_4$.

It follows that 
\begin{equation}\label{minor2}{\det}_q=\sum_{M_4\in{\mathcal
M}_4}\sum_{G\in{\mathcal
R}_{1,2,3}}P_{G,M_4}Z^GM_4,
\end{equation}where each $P_{G,M_4}\in{\mathbb Z}[q^2]$. Then, because 
reordering elements from ${\mathcal R}_{1,2,3}$ does not introduce terms from
${\mathcal R}_4$,
\begin{equation}\label{minor3}Z(A_0){\det}_q=\sum_{M_4\in{\mathcal
M}_4}\sum_{H\in{\mathcal
R}_{1,2,3}}\hat{P}_{H,M_4}Z^HM_4
\end{equation}for some elements $\hat{P}_{H,M_4}\in{\mathbb Z}[q^2,q^{-2}]$. At
the same time, by the induction hypothesis,
\begin{equation}\label{minor44}Z(A_0){\det}_q=\sum_{L\in {\mathcal R}_4
}\sum_{K\in{\mathcal R}_{1,2,3}}\tilde{c}_{L,K}Z(K+L),
\end{equation}
where $\tilde{c}_{L,K}=1$ for the unique configuration corresponding to
$Z(A_0+I)$ and in all other cases, if non-zero,  $\tilde{c}_{L,K}=\pm
q^{2\gamma_{K,L}}$ where $\gamma_{K,L}\in {\mathbb N}$. Here, each $K+L$ has the same row and column sums as $A_0+I$, and $K+L\leq A_0+I$. The expression $Z^HM_4$
in (\ref{minor3}) is a sum of monomials $\pm q^{2p_i}Z^HZ^{S_{4,i}}$ corresponding
to $M_4=\sum_i\pm q^{2p_i}Z^{S_{4,i}}$ as a quantum minor. Furthermore, $\forall i: p_i\in{\mathbb Z}_+$ and
$p_i\in{\mathbb N}$ for all but one $i$. Due to the configurations, the normalization factors
$N(H+S_{4,i})$ are easily seen to be  independent of $i$ and may thus for instance be computed for
the unique $i$ for which $p_i=0$. Thus, $N(H+S_{4,i})=N(H)q^{-\underline{ro}(H;M_4)-\underline{co}(H;M_4)}$. The symbols
$\underline{ro}(H;M_4)$ and $\underline{co}(H;M_4)$  denote the row sums,
respectively column sums, of $H$ corresponding to the rows, respectively
columns, of $M_4$. Thus, for each $M_4$ in (\ref{minor3}) we get a sum of the form $\pm\hat P_{H,M_4}N(H)^{-1}q^{\underline{ro}(H;M_4)+\underline{co}(H;M_4)}q^{2p_i}Z(H+S_{4,i})$. Each of these terms must correspond to a term in (\ref{minor44}) where we have full information about the positivity of the powers of $q$. Notice that $\pm\hat P_{H,M_4}$ is independent of $i$. Indeed, $\hat P_{H,M_4}=c_{H,M_4}
q^{2p_{H,M_4}}q^{-\underline{ro}(H;M_4)-\underline{co}(H;M_4)}N(H)$  for some constant $c_{H,M_4}\in\{-1,0,1\}$ and some element
${p}_{H,M_4}\in{\mathbb Z}_+$. It follows from this that we have a formula

\begin{equation}\label{minor5}Z(A_0){\det}_q=\sum_{M_4\in{\mathcal
M}_4}\sum_{H\in{\mathcal
R}_{1,2,3}}c_{H,M_4}
q^{2p_{H,M_4}}q^{-\underline{ro}(H;M_4)-\underline{co}(H;M_4)}Z(H)M_4.
\end{equation} Each summand in
$q^{-\underline{ro}(H;M_4)-\underline{co}(H;M_4)}Z(H)M_4$ is normalized. With the
exception of one pair $(H_s,M_{4,s})$ where
${p}_{H_s,M_{4,s}}=0$, we have furthermore that ${p}_{H,M_4}\in {\mathbb N}$
when $c_{H,M_4}\neq0$.

Observing  that  ${\det}_q$ is central, we can insert it in
any position we prefer. Returning to (\ref{reorder}) we will therefore consider
$Z(A_0)\cdot det_q\cdot 
q^{-\alpha}Z_{r+1,c}^{a_{r+1,c}}$. In view of (\ref{minor5}) we need to focus on
the
rewriting of expressions of the form $ M_4Z_{r+1,c}^{a_{r+1,c}}$ and, in
particular, to carefully keep 
track of the $q$ factors we pick up. This is the only place where negative
exponents might originate. 

Four different situations may occur:

\begin{itemize}
\item{1)} $M_4$ has all row numbers greater than $r+1$ and all column numbers
greater than $c$.
\item{2)} $M_4$ has all row numbers greater than $r+1$ but a  column number
equal to $c$.
\item{3)} $M_4$ has a row number equal to $r+1$ but all column numbers greater
than $c$.
\item{4)} $Z_{r+1,c}$ occurs in $M_4$ {{}( and hence commutes with $M_4$.)}
\end{itemize}

These cases can be dealt with, and this is, indeed, the reason for the chosen decomposition: In all cases, what has to be considered are the terms 
\begin{equation}
(\star) = c_{H,M_4}q^{2p_{H,M_4}}q^{-\underline{ro}(H;M_4)-\underline{co}(H;M_4)}Z(H)M_4q^{-\alpha}Z^a_{r+1,c}.
\end{equation}

In cases 2, 3), and 4), $M_4$ quasi-commutes with $Z^a_{r+1,c}$, and $(\star)$ becomes one of the following:

\begin{eqnarray*}
(\star)_2=(\star)_3&=& c_{H,M_4}q^{2p_{H,M_4}}q^{-\underline{ro}(H;M_4)-\underline{co}(H;M_4)}q^{-\alpha}q^{-2a}Z(H)Z^a_{r+1,c}M_4\\
(\star)_4&=& c_{H,M_4}q^{2p_{H,M_4}}q^{-\underline{ro}(H;M_4)-\underline{co}(H;M_4)}q^{-\alpha}Z(H)Z^a_{r+1,c}M_4.
\end{eqnarray*}

The quantum minor $M_4$ is a linear combination of monomials with coefficients which are non-negative powers of $(-q^2)$. These powers are no problem (they just become part of the coefficients $q^{2\gamma_B}c_B$ in (\ref{det-for})), so what  need to be dealt with are the expressions resulting from replacing $M_4$ in $(\star)_2$,   $(\star)_3$, and  $(\star)_4$ by any of the monomials $A$ from $M_4$.

As already noted, the non-zero positions in $Z^a_{r+1,c}A$ are ordered correctly according to the lexicographic ordering, so $Z^a_{r+1,c}A=Z^B$ for some $B\in M_n({\mathbb Z}_+)$.

In case 2), $N(B)=q^{-a}$, so $Z^a_{r+1,c}A=q^aZ(B)$ and the term from $(\star)_2$ we need to analyze is
\begin{equation}
(\star\star)_2= c_{H,M_4}q^{2p_{H,M_4}}q^{-\underline{ro}(H;M_4)-\underline{co}(H;M_4)}q^{-\alpha}q^{-a}Z(H)Z(B).
\end{equation}

There are no positions where $H$ and $B$ both have non-zero entries. Moreover, if $H$ and $B$ have non-zero entries in the respective positions $(i,j)$ and $(k,l)$ with  $(k,l)>(i,j)$, then the corresponding powers $Z^h_{ij}$ in $Z^H$ and $Z^b_{kl}$ in $Z^B$ commute. This means $Z^HZ^B=Z^C$, where $C=H+B$. Observe that 
\begin{equation}
N(C)=N(H)N(B)q^{-\underline{ro}(H;M_4)-\underline{co}(H;M_4)}q^{-a\underline{ro}_{r+1}(H)-a\underline{co}_c(H)}.
\end{equation}
The matrix $H$ consists of $A_0$ plus the complementary part of a permutation matrix whose other part is a permutation matrix corresponding to a term in $M_4$. So,  $\underline{ro}_i(H)=\underline{ro}_i(A_0)$ for rows $i$ that occur in $M_4$ and $\underline{ro}_i(H)=\underline{ro}_i(A_0)+1$ otherwise. Similarly for the columns. In case 2), this means $\underline{ro}_{r+1}(H)=\underline{ro}_{r+1}(A_0)+1$ and $\underline{co}_c(H)=\underline{co}_c(A_0)$. Now, $a\underline{ro}_{r+1}(H)+a\underline{co}_c(H)=\alpha+a$, so

\begin{equation}
N(C)=N(H)N(B)q^{-\alpha-a-\underline{ro}(H;M_4)-\underline{co}(H;M_4)}
\end{equation}
and $Z(C)=q^{-\alpha-a-\underline{ro}(H;M_4)-\underline{co}(H;M_4)}Z(H)Z(B)$. Therefore
\begin{equation}
(\star\star)_2=c_{H,M_4}q^{p_{H,M_4}}Z(C),
\end{equation}
the correct form for the right hand side of (\ref{det-for}).

Case 3) is symmetric.

In case 4), $N(B)$ is either $1$ or $q^{-2a}$, depending on whether or not $A$ has a $1$ in position $(r+1,c)$. This time, $\underline{ro}_{r+1}(H)=\underline{ro}_{r+1}(A_0)$ and $\underline{co}_c(H)=\underline{co}_c(A_0)$, so $a\underline{ro}_{r+1}(H)+a\underline{co}_c(H)=\alpha$ and we end up with, setting $C=H+B$,
\begin{equation}
(\star\star)_4=c_{H,M_4}q^{2p_{H,M_4}}\cdot[1\textrm{ or }q^{4a}]\cdot Z(C),
\end{equation}
in the correct form.

In case 1) we have a reinterpretation of (\ref{commu0}):
\begin{equation}\label{minor4}[Z_{r+1,c}^a,M_4]=-q^2(1-q^{-4a})Z_{r+1,c}^{a-1}T,
\end{equation} 
where $Z_{r+1,c}T=q^4TZ_{r+1,c}$, and $a=a_{r+1,c}$. The factor $T$ will be discussed shortly.

Notice that the left hand side of (\ref{minor4}) evidently is skew under the bar
operator. Thus, it follows  that
$\overline{Z_{r+1,c}^{a-1}T}=q^{-4a+4}Z_{r+1,c}^{a-1}T$. We have $T=q^{-2}N_1$ in terms of (\ref{commu0}), so each monomial
$Z^{Y_i}=Z(Y_i)$ in $T$ occurs with a factor $\pm q^{2p_i}$ with $p_i\in{\mathbb
Z}_+$. However, we must utilize even finer details of $T$. Specifically, we
may assume that the monomial summands of $T$ each have a contribution $Z_{x,c}$
with $x>r+1$ and a contribution $Z_{r+1,y}$ with $y>c$. Furthermore, $T$ is
ordered according to the lexicographic ordering column-row and a  factor of
$q^2$ is taken out of the original determinantal expression which involves
expressions  $(-q^2)^{\ell}$, where  $\ell\geq1$. It is clearly the term with
$q^{2-4a}$ we must be able to handle. Before addressing this, we remark that the
term $Z_{r+1,c}^aM_4$ from the commutator is handled by the same argument as in
cases 2), 3), and 4).

We know from the construction that each $K+L$ in (\ref{minor44}), appearing
with a non-zero coefficient, compared to $A_0$ has an additional element in
each row and column coming from the various summands in the determinant. With
the given $M_4$ we then know that the extra element $W_{r+1,u}$ in the $(r+1)$th
row must have $u<c$ and
the extra element $W_{v,c}$ in the $c$th column must have $v<r+1$.

The above observations easily imply that \begin{eqnarray}\label{baragain}
\overline{Z(H)Z_{r+1,c}^{a-1}T}&=&q^{-4a+4}Z_{r+1,c}^{a-1}T\overline{Z(H)}
\nonumber\\\label{long}
&=&q^{-4a+4}q^{-2\alpha}q^{-4a}q^{-2\underline{ro}(H;M_4)-2\underline{co}(H;M_4)}
Z(H)Z_{r+1,c}^{a-1}T \\&+&\textrm{ lower order terms}.\nonumber
\end{eqnarray}

The term we have to control is

$$X=q^2q^{-4a}q^{-\alpha}q^{-\underline{ro}(H;M_4)-\underline{co}(H;M_4)}Z(H)Z_{
r+1,c}^{a-1}T.$$

Equation (\ref{long}) implies that the term in $X$ coming from the leading term
in $T$ is normalized. The other terms are then positive powers of $q^2$ times
normalized elements {{}as follows by arguments similar to those for the cases 2), 3), and 4).}

This completes the proof. \qed

\bigskip

\section{Covariant Minors and the dual canonical basis}

Let us consider an $n\times n$ matrix $X\in M_n({\mathbb Z}_+)$  decomposed into
\begin{equation}
  \label{eq:dec}
  X=
  \left(\begin{array}{cc}A&B\\C&D    
  \end{array}\right).
\end{equation}
We assume furthermore that $C$ is square of size $s$. We denote
the  $s\times s$ quantum minor corresponding to the lower left corner
by $I_{s,ll}$.
\begin{Lem}\label{lem}
Let $b(X)$ be an element of the dual canonical basis with $X$ given as
in (\ref{eq:dec}). Then
\begin{equation}
  \label{eq:basis}
  b(X)I_{s,ll} = q^{(S(A)-S(D))}b(\tilde X),
\end{equation}
with $\tilde X=\left(\begin{array}{cc}A&B\\C+I_s&D    
  \end{array}\right)$. Here, $I_s$ is the $s\times s$ identity matrix,
while $S(A)$ and $S(D)$ denote the sum of all entries in $A$ and $D$,
respectively.  
\end{Lem}

\pof It is easy to see that $Z^X=Z^AZ^BZ^CZ^D$. Suppose then, by
Proposition~\ref{basis}, that
\begin{equation}
  \label{eq:x-dec}
b(X)=Z(X)+\sum_{X^\prime< X}c_{X^\prime}(X)Z(X^\prime),
\end{equation}
with $c_{X^\prime}(X)\in q^2{\mathbb Z}[q^2]$. Furthermore, each $X^\prime$ has the same row and column sums as $X$. Here we set
$X^\prime=
  \left(\begin{array}{cc}A^\prime&B^\prime\\C^\prime&D^\prime    
  \end{array}\right)$  and then  
$$Z(X^\prime)=C_{A^\prime,B^\prime,C^\prime,D^\prime}Z^{A^\prime}Z^{B^\prime}Z^{
C^\prime}Z^{
D^\prime}$$where
\begin{eqnarray*}&C_{A^\prime,B^\prime,C^\prime,D^\prime}=\\&
  N(A^\prime)N(B^\prime)N(C^\prime)N(D^\prime)N(A^\prime,B^\prime)
  N(A^\prime,C^\prime) N(B^\prime,D^\prime)N(C^\prime,D^\prime).\end{eqnarray*}

The factors $N(A^\prime)$, $N(B^\prime)$,
$N(C^\prime)$, and $N(D^\prime)$ are given by (\ref{norma}) as are the  factors of mixed summands;
$N(A^\prime,B^\prime)=q^{-
  \sum_j \underline{ro}_j(A^\prime)\underline{ro}_j(B^\prime)}, N(A^\prime,C^\prime)=q^{-
  \sum_j \underline{co}_j(A^\prime)\underline{co}_j(C^\prime)}$, $N(B^\prime,D^\prime)$, and
$N(C^\prime,D^\prime)$ (there will be no non-trivial factors
$N(A^\prime,D^\prime)$ and $N(B^\prime,C^\prime)$).

Looking at row sums we have $S(A^\prime)+S(B^\prime)=S(A)+S(B)$, and looking at
column sums we have $S(B^\prime)+S(D^\prime)=S(B)+S(D)$. Thus, for the
matrices in the right hand side of (\ref{eq:x-dec}) we have
$-S(A^\prime)+S(D^\prime)=-S(A)+S(D)$. Then,
\begin{eqnarray}
  \label{basis2}
  &q^{(-S(A)+S(D))}
  b(X)I_{s,ll}\\&=\sum_{A^\prime,B^\prime,C^\prime,D^\prime}
  q^{(-S(A^\prime)+S(D^\prime))} 
c_{X^\prime}(X)C_{A^\prime,B^\prime,C^\prime,D^\prime}
Z^{A^\prime}Z^{B^\prime}Z^{C^\prime}Z^{D^\prime}I_{s,ll}.                 \nonumber
\end{eqnarray} 
It follows from Proposition~\ref{cminor} (a transposed version thereof) that
$$I_{s,ll}Z^{A^\prime}Z^{B^\prime}Z^{C^\prime}Z^{D^\prime}=q^{-2S(A)+2S(D)}Z^{
A^\prime}Z^{B^\prime}Z^{C^\prime}Z^{D^\prime}I_{s,ll}.$$ Thus, the left hand
side, and hence both sides of
(\ref{basis2}) are
bar invariant. Now consider a term in (\ref{basis2}) of the form
\begin{eqnarray*}q^{(-S(A^\prime)+S(D^\prime))}
C_{A^\prime,B^\prime,C^\prime,D^\prime}
Z^{A^\prime}Z^{B^\prime}Z^{C^\prime}Z^{D^\prime}I_{s,ll}\\=q^{(-S(A^\prime)-S(D^\prime))}
C_{A^\prime,B^\prime,C^\prime,D^\prime}
Z^{A^\prime}Z^{B^\prime}Z^{C^\prime}I_{s,ll}Z^{D^\prime}.\end{eqnarray*}

Here, $Z(C^\prime)I_{s,ll}=N(C^\prime)Z^{C^\prime}I_{s,ll}$. By
Theorem~\ref{thm} this equals $Z(C^\prime+I_s) +
\sum_{C''<C^\prime +I_s}f_{C''}Z(C'')$ and for all $C''$, $f_{C''}$ a polynomial
in
$q^2{\mathbb Z}[q^2]$.
Notice that each $C''$ has the same row and column sums as
$C^\prime+I_s$ and that  $\forall i:
\underline{ro}_i(C^\prime+I_s)=\underline{ro}_i(C^\prime)+1$ and, similarly,
$\forall j:
\underline{co}_j(C^\prime+I_s)=\underline{co}_j(C^\prime)+1$. But then $
q^{-S(A^\prime)-S(D^\prime)}N(A^\prime,C^\prime)N(C^\prime,
D^\prime)=N(A^\prime,C'')N(C'',D)$. Thus, the right hand
side of (\ref{basis2}) is a sum of terms $g_{i}Z(Y_i)$ with
$Y_i=\left(\begin{array}{cc}A^\prime&B^\prime\\C^{\prime\prime}&D^\prime    
  \end{array}\right)$ and $g_i\in {\mathbb Z}[q^2]$. Precisely
the term
with  $A^\prime=A,\ B^\prime=B, \ C^\prime=C$, and $D^\prime=D$ has a
coefficient
$g_i=1$, all other coefficients are in $q^2{\mathbb Z}[q^2]$. Thus the right
hand side has the right expansion
properties, hence is a member of the dual canonical basis corresponding to the
stated element $\tilde X$. \qed

\medskip


Let us instead consider an $n\times n$ matrix $X\in M_n({\mathbb Z}_+)$ 
decomposed into
\begin{equation}
  \label{eq:dec-0}
  X=
  \left(\begin{array}{cc}0&B\\C&D    
  \end{array}\right),
\end{equation}
where we now  assume  that $D$ is square of size $s$. We denote
the  $s\times s$ quantum minor corresponding to the lower right hand
corner (as occupied by $D$)
by $I_{s,lr}$.
\begin{Lem}\label{samelem}
Let $b(X)$ be an element of the dual canonical basis with $X$ given as
in (\ref{eq:dec-0}). Then
\begin{equation}
  \label{eq:basis-now}
  b(X)I_{s,lr} = q^{(S(B)+S(C))}b(\tilde X),
\end{equation}
with $\tilde X=\left(\begin{array}{cc}0&B\\C&D+I_s    
  \end{array}\right)$. As before, $I_s$ is the $s\times s$ identity matrix,
while $S(B)$ and $S(C)$ denote the sum of all entries in $B$ and $C$,
respectively.   
\end{Lem}

The proof follows the same lines as that of Lemma~\ref{lem} and is
omitted. By Proposition~\ref{simple} the following case encompasses the two
former. The proof is omitted for similar reasons.

\begin{Lem}\label{samelem2}
Let $b(X)$ be an element of the dual canonical basis with $$X=\left(
  \begin{array}{ccc}0&B_1&B_2\\C_1&D&C_2\\G_1&G_2&0
  \end{array}\right)$$where $D$ is $s\times s$. Then, if $I_{s,cc}$ denotes the
$s\times s$ quantum minor corresponding to the position of $D$,
\begin{equation}
  \label{eq:basis-now2}
  b(X)I_{s,cc} = q^{S(B_1)+S(C_1)-S(C_2)-S(G_2)}b(\tilde X),
\end{equation}
with $\tilde X=\left(
  \begin{array}{ccc}0&B_1&B_2\\C_1&D+I_s&C_2\\G_1&G_2&0
  \end{array}\right)$. As before, $I_s$ is the $s\times s$ identity matrix.
\end{Lem}
\label{bas}
\medskip
\begin{Rem}\label{same3}Using Proposition~\ref{simple} it follows that analogous
results hold for  the configurations 
$$X=\left(
  \begin{array}{ccc}C_1&D&C_2\\G_1&G_2&0
  \end{array}\right) \ ,\ \tilde X=\left(
  \begin{array}{ccc}C_1&D+I_s&C_2\\G_1&G_2&0
  \end{array}\right)$$
  and 
  $$X=\left(
  \begin{array}{cc}0&B_2\\D&C_2\\G_1&G_2
  \end{array}\right)\ ,\ \tilde X=\left(
  \begin{array}{cc}0&B_2\\D+I_s&C_2\\G_1&G_2
  \end{array}\right),$$ 
in which the matrix $X$ is not necessarily square. Similarly, the transposed
cases,  where the $0$ matrix is in the opposite corner, are covered.  
\end{Rem}

\medskip

\section{Broken line constructions}\label{6}

Consider the $m\times n$ quantum matrix algebra ${\mathcal
  O}_q(M(m,n))$.  In this section, all elements $Z_{i,j}$ and all quantum
minors are elements of this algebra.  

\begin{Def} {\bf A
broken line in } $M_{m,n}({\mathbb Z}_+)$ is a path in ${\mathbb
N}\times {\mathbb N}$ starting at $(1,n)$  and terminating at $(m,1)$. We will
occasionally also refer to this as a broken line from $(1,n)$ to $(m,1)$. It
must satisfy furthermore that it alternates between horizontal and vertical
segments while passing through
smaller column numbers (in the horizontal
direction)  and bigger row numbers (in the vertical direction). 

Unless we are in the
extreme cases
$(1,n)\mapsto (1,1)\mapsto(m,1)$ or $(1,n)\mapsto (m,n)\mapsto(m,1)$, this will
divide the indices
$(i,j)$ into 3 disjoint sets $S_L, L$, and $T_L$. Here, $S_L$ is the set of
points above the line (when there are 3 subsets, we have that $(1,1)$ is
above the
line), $L$ is the line itself, and $T_L$ is the set of  points
below the line. 
\end{Def}

\begin{Rem}\label{vein}
A broken line is determined by a double partition 
$$1=i_1\leq i_2\leq i_3\leq \dots \leq i_s=m\ \textrm{ and } n=j_1\geq j_2\geq
j_3\geq\dots\geq j_s=1,$$
such that the corners in the line $L$ are $(i_t,j_t); t=1,2,\dots, s$. This,
naturally, dictates that in the partitions, precisely every second inequality is
sharp. Furthermore, if in a given position, one is sharp, then the other is not,
and vice versa. 

In a similar vein, the broken line is given by a double flag variety. 
\end{Rem}

For a given broken line $L$, we  now construct a family ${\mathcal V}_L$
with $mn$ elements
consisting of certain quantum minors: (It will be proved below that all members 
$q$-commute.) For points in $(i,j)\in T_L\cup L$ we take the biggest
quantum minor
having
its bottom right corner in $(i,j)$ and completely contained in  $T_L\cup
L$. One can also say that it is the biggest quantum minor consisting of
adjacent rows and columns (we call such a quantum minor {\bf solid}) and which
contains
$(i,j)$ as well as
points from $L$ but no points from $S_L$. The line $L$ is thus represented by
points, that is, $1\times1$ matrices. For the points in $S_L$ we do something
else: For $(i,j)\in 
S_L$ we
take the biggest quantum minor consisting of adjacent
 rows and columns and which
contains $(i,j)$ in the upper left corner (all other rows have
numbers bigger than $i$ and all columns have numbers bigger than $j$). Notice
that with $L$ fixed, each quantum minor in ${\mathcal V}_L$ corresponds uniquely
to a point
$(i,j)$. By
{\em the quantum minor corresponding to a point} we then mean this quantum
minor.

\medskip

The first important observation is:
\begin{Prop}\label{a-l} Any quantum minor corresponding to a point in $S_L$
$q$-commutes
  with any $Z_{i,j}$ for which $(i,j)\notin S_L$.
\end{Prop}

This follows immediately from Proposition~\ref{cminor}. 

{{}
\begin{Prop}\label{q-cor} Let $M=M_{a,b}(k)$ be a $k\times k$ quantum minor
with
  upper left corner in $(a,b)$ and lower right hand corner in
  $(a+k-1,b+k-1)$, and such that $M$ is inside the $m\times n$ quantum
  matrices. Refer to the 9 different positions of a pair $(i,j)$  relative
to $M$ as $NW(M), N(M), NE(M)$, $ \dots, I(M),\dots$, $SE(M)$ such that $NW(M)$ is ($a>i$ and $b>j$) and SE(M) is
($i>a+k-1$ and $j>b+k-1$). Here $I(M)$ denotes the position of the indices of $M$.
Then $Z_{i,j}$
$q$-commutes with $M$ {\em unless} $(i,j)$ is in  $NW(M)$ or in $SE(M)$. For the remaining
pairs, in the $q$-commutation formulas  $Z_{i,j}M=q^{2p_{i,j}}MZ_{i,j}$,
$p_{i,j}$
depends only on the relative positions. Indeed, $p_{i,j}=1$ for $(i,j)$ in $W(M)\cup N(M)$, $p_{i,j}=0$ for $(i,j)$ in $I(M)\cup SW(M)\cup NE(M)$,  and $p_{i,j}=-1$ for
$(i,j)$ in $S(M)\cup E(M)$.
  
\end{Prop}

\proof  With the exception of NW(M) and SE(M), the
$q$-commutation relation may be seen as taking place inside a
smaller matrix algebra in which $M$ is a covariant quantum minor. \qed

\begin{Prop} All members of ${\mathcal V}_L$ $q$-commute. \end{Prop}

\proof Let $A,B\in {\mathcal V}_L$. Let ${\mathcal E}(A)$ denote the entries of $A$. It follows  by inspection that, after possibly interchanging $A$ and $B$, only the following cases may occur: ${\mathcal E}(A)\subseteq S(B)\cup SW(B)\cup W(B)\cup I(B)$,    ${\mathcal E}(A)\subseteq N(B)\cup NE(B)\cup E(B)\cup I(B)$, or ${\mathcal E}(B)\subseteq {\mathcal E}(A)$. Some situations involving fewer sets like ${\mathcal E}(A)\subseteq S(B)\cup I(B)$ or ${\mathcal E}(A)\subseteq S(B)$ may also occur, while others may be prohibited  due to the configuration at hand.  All non-trivial cases are treated in the same way and it suffices to consider the
very first of these. Let $B$
be fixed and consider the expansion of $A$ into a linear combination of
monomials of the form $Z_{i_1+1,j_1+\sigma(1)}\cdots 
Z_{i_1+r,j_1+\sigma(r)}$ for
some $\sigma\in S_r$. Using Proposition~\ref{a-l} and
Proposition~\ref{q-cor} we obtain the following: If  $W_\sigma(B)$ and 
$S_\sigma(B)$ denote the number of terms $Z_{i_1+i,\sigma(i_1)+i}$ to the west,
respectively to the south, of $B$, the given monomial will $q$-commute with $B$
with a factor $q^{2(W_\sigma(B)-S_\sigma(B))}$. It is easily seen that
$W_\sigma(B)-S_\sigma(B)$ is independent of $\sigma$, and thus the claim
follows. \qed}

This result also follows from \cite[Theorem~1]{sco}.

\begin{Rem}One gets a similar family by interchanging $S_L$ and $T_L$.
Colloquially speaking, if one allows $L$ to vary, this can be accomplished by a
reflection mapping $m\times n$ matrices to $n\times m$ matrices while
interchanging rows and columns.
\end{Rem}

\medskip

\begin{Def} We introduce a partial ordering of the broken lines:
$$L_1\leq L_2\Leftrightarrow S_{L_2}\subseteq S_{L_1}.$$
\end{Def}

In this ordering, the line $L^+$ corresponding to the empty set:
$(1,n)\rightarrow (1,1)\rightarrow(m,1)$, is the unique
maximal element, and the line $L^-$ corresponding to $T_L=\emptyset$:
$(1,n)\rightarrow (m,n)\rightarrow(m,1)$, is the unique
smallest element.

In the extreme case of $L^+$, the $q$-commuting quantum minors in the
corresponding family are the
following:

\begin{enumerate}

\item For $i\ge j$, $\xi^{\{i-j+1, i-j+2,\dots,i\}}_{\{1,2,\dots,j\}}$,
\newline
\item For $j>i$, $\xi^{\{1,2,\dots,i\}}_{\{j-i+1, j-i+2,\dots,j\}}$.
\end{enumerate}

\medskip

In the sequel, we shall consider the
following more general family ${\mathcal
V}_{L_1,L_2}\subset {\mathcal V}_{L_1}$ where, clearly, ${\mathcal
V}_L={\mathcal V}_{L,L^+}$:

\begin{Def}Let $L_1,L_2$ be broken lines with $L_1\leq L_2$. The family  ${\mathcal
V}_{L_1,L_2}$ is the subfamily of ${\mathcal V}_{L_1}$ that    corresponds to
the points in $T_{L_2}\cup L_2$.
\end{Def}

\medskip

\begin{Def}Given a broken line $L$, let ${\mathcal
    O}_q(T_L\cup L)$ denote the sub-algebra of ${\mathcal O}_q(M(m,n))$
  generated by the $Z_{i,j}$ for which $(i,j)\in T_L\cup L$. We will refer to
the members of ${\mathcal V}_L$ as {\em variables}. 
Let ${\mathcal  V}_L^-$ denote the set of variables in ${\mathcal V}_L$
corresponding to the points in $L\cup T_L$, and let ${\mathcal C}^-_{L}$ denote
the set of variables  in ${\mathcal V}_L$ for the points in $L^-$ $(\subseteq
(T_L\cup L))$.
Analogously, let ${\mathcal V}_L^+$
denote the set of 
  variables in  ${\mathcal V}_L$ corresponding to the points in $S_L$.
 
In the
following we will consider cluster algebra constructions inside an ambient
space which is either i)  the skew field of fractions ${\mathcal
F}_L$ constructed from ${\mathcal O}_q(M(m,n))$  (or, equivalently, 
${\mathcal V}_L$) and where ${\mathcal V}_L$ is part of an initial seed or ii) 
the skew field of fractions ${\mathcal
F}^-_L$ constructed from ${\mathcal
    O}_q(T_L\cup L)$ (or, equivalently, ${\mathcal V}^-_L$) and where
${\mathcal V}^-_L$ is part of an initial seed. 
\end{Def}

Proposition~\ref{a-l} can be stated as the fact that any
variable in ${\mathcal V}_L^+$ is covariant with respect to the full
sub-algebra ${\mathcal O}_q(T_L\cup L)$. The algebra ${\mathcal O}_q(T_L\cup L)$
has previously been studied in \cite[section 3]{llr}.

\medskip

\medskip

\begin{Thm}Let ${\mathcal V}_L$ denote the family of $mn$ $q$-commuting
  quantum minors constructed from a broken line as above. Then  up to
  multiplication by a power of $q$,  any monomial in
  the members of  ${\mathcal V}_L$ is a member of the dual canonical
  basis.
\end{Thm}

\pof The main tool is Lemma~\ref{samelem2}, but Proposition~\ref{simple} is also
important, c.f. the remark following Lemma~\ref{samelem2}. Consider then a
monomial. Rewrite it, if necessary, in such a way that the factors coming from
the points on $L$ are furthest to the left. Then place to the right of these the
$2\times2$-minors corresponding to the points one step below the line.
Continue in this way  until all the factors corresponding to the points on, or
below, the line are positioned. While continuing to add from the right, 
order the factors coming from $S_L$ in a similar fashion and such that the
factor corresponding to the position $(1,1)$ is furthest to the right. The finer
order is not important. We view the monomial as the result of a sequence of
multiplications from the right by minors according to this ordering.
Inductively, we may at each step $r$ in the sequence assume that what we are
multiplying the minor onto is some $q^{2p_r}b(X_r)$. The start
is clearly
trivial.  When we add minors below
the line, $X_r$ is all the time of the form as given in
Lemma~\ref{samelem2}. After that the
zero  in the lower right corner disappears and  we apply
Remark~\ref{same3} instead. The
result follows. \qed

\bigskip

\begin{Def}\label{closest}For a given line $L$, we say that the line $L_1$ is a
  closest bigger line to $L$ if $L<L_1$ and there is no other line $L_2$ such
  that $L<L_2<L_1$. In this case, if $L=(1,n)\rightarrow\dots\rightarrow
(f,d)\rightarrow(c,d)\rightarrow(c, g)\rightarrow\dots\rightarrow(m,1)$, then
$L_1=(1,n)\rightarrow\dots\rightarrow
(f,d)\rightarrow(c-1,d)\rightarrow(c-1,d-1)\rightarrow(c,d-1)\rightarrow(c,
g)\rightarrow\dots\rightarrow(m,1)$ for some such ``corner'' $(f,d)\rightarrow
  (c,d)\rightarrow (c,g)$, where we, naturally, also allow $f=c-1$
  and $g=d-1$.
  
  We will call the given corner of $L$ {\bf convex} and the resulting corner of
$L_1$ {\bf concave}. We will also write
$$L_1=L\uparrow (c,d) \textrm{ or }L=L_1\downarrow (c-1,d-1).$$ 
\end{Def}

\bigskip

\subsection{Key technical results}

Focus on a position $(i_0,j_0)$ inside the quantum  matrix algebra ${\mathcal O}_q(M(n_0,r_0))$.  
Consider the sub-algebra $M={\mathcal O}_q^{i_0,j_0}(M(s))$
generated by the variables $Z_{i_0+a,j_0+b}$ with $0\leq a,b \leq
s-1$ where $s$ is the biggest positive integer such that $Z_{i_0+s-1,j_0+s-1}\in
{\mathcal O}_q(M(n_0,  r_0))$. Naturally, this sub-algebra is  isomorphic to
${\mathcal
O}_q(M(s))$ in which we number the rows and columns as $0,1,\dots, s-1$.
Assume $s\geq2$. Inside
$M$ are the quantum minors
$Y_r=Y_r^{(s-2)}=\xi_{\{1,\dots,s-1\}}^{\{0,1,\dots,s-2\}}$,
$Y_l=Y_l^{(s-2)}=\xi_{\{0,\dots,s-2\}}^{\{1,\dots,s-1\}}$,
$X_o=X_o^{(s-2)}=\xi_{\{1,\dots,s-2\}}^{\{1,\dots,s-2\}}$,
$X_t=X_t^{(s-2)}=\xi_{\{0,\dots,s-2\}}^{\{0,\dots,s-2\}}$,
$X_b=X_b^{(s-2)}=\xi_{\{1,\dots,s-1\}}^{\{1,\cdots,s-1\}}$, and 
$D=D^{(s-2)}=\xi_{\{0,1,\dots,s-1\}}^{\{0,1,\dots,s-1\}}$. The last is 
just the full quantum determinant in $M$. $X_o^{(0)}$ is defined as the constant $1$.

\medskip

\begin{Def}
A set $\{X_t,X_b, D, X_o,Y_l,Y_r\}\subset {\mathcal O}_q(M(n_0, 
r_0))$ whose elements are given as above for some $i_0,j_0,s\in {\mathbb N}$ will be called an ${\mathcal M}$-set.
\end{Def}

We have the following facts which follow by direct computation:

{{}\begin{Lem}The elements $Y_\ell$ and $Y_r$ commute. The elements $D,X_o$, and $X_b$ commute, and we have\label{red1}
\begin{equation}\overline{X_oDX_b^{-1}}=X_oDX_b^{-1}\textrm{ and }\ \overline{q^2Y_rY_\ell X_b^{-1}}=q^2Y_rY_\ell X_b^{-1}.
\end{equation}

\end{Lem}

\begin{Cor}\label{level}
The element $X_o^{-1}D^{-1}Y_rY_l$ commutes with all elements
$Z_{i,j}\in M$ with the exception of $Z_{i_0,j_0}$
and $Z_{i_0+s-1,j_0+s-1}$. In particular, it commutes with the
quantum minors $X_o^{(a)}=\xi_{\{1,\dots,a\}}^{\{1,\cdots,a\}}$ for $a=1,\dots,
s-2$.
It
$q$-commutes with the quantum minors
$X_b$, and $X_t$
according to
 \begin{eqnarray}X_bX_o^{-1}D^{-1}Y_rY_l&=&q^{-4}X_o^{-1}D^{-1}
Y_rY_lX_b\textrm{ and }\\X_tX_o^{-1}D^{-1}Y_rY_l&=&q^{4}X_o^{-1}D^{-1}
Y_rY_lX_t.\nonumber\end{eqnarray}
\label{4.1} 
\end{Cor}
}

\bigskip

The following result will allow us to construct the $B$ matrices of
the compatible pairs. It follows by easy computation.

\begin{Cor}\label{import}
Let $\Lambda_{\small}$ be defined as the $\Lambda$-matrix of Laurent quasi-polynomial algebra generated by the
variables $X_b,X_o,D,Y_r,Y_l$ and their inverses. Then, 
$$\Lambda \left(\begin{array}{c}0\\-1\\-1\\1\\1\end{array}\right)=
\left(\begin{array}{c}-4\\0\\0\\0\\0\end{array}\right).$$ 
\end{Cor}

\medskip
 
The following was proved by
Parshall and Wang in \cite[Theorem~5.2.1]{pw} but is also a special
case of \cite{good} Theorem~6.2. 

\begin{Prop}[Parshall and Wang] \label{4.3}
$$X_tX_b-X_bX_t= (q^2-q^{-2}) Y_rY_l.$$
\end{Prop}

We wish to strengthen this result considerably, namely to the following equation
which will play an important
role later when we consider the quantum mutations in certain
directions. 

\begin{Thm} \label{4.4}
$$X_tX_b= X_oD + q^2Y_rY_l.$$
\end{Thm}

\pof We first observe that by Theorem~\ref{thm} and Lemma~\ref{lem},  $X_oD$
and
$Y_rY_l$ are members of the dual canonical basis;  $X_oD=b(A_1)$
and $Y_rY_l=b(A_2)$ for some specific matrices $A_1,A_2$ with $A_2<A_1$. Notice that when cut down to the $s\times s$ block where it is located, 
\begin{equation}A_1=\left(\begin{array}{ccccc}1&&\cdots&&0\\&2\\\vdots&&\ddots&&\vdots\\&&&2\\0&&\cdots&&1\end{array}\right).\end{equation}We
consider the expansion of $X_tX_b$ onto the dual canonical basis;
\begin{equation}
  X_tX_b=\sum_i c_i(q)b(C_i).
\end{equation}
The coefficients $c_i$ are in ${\mathbb Z}[q^{2},q^{-2}]$. Since $X_tX_b=Z^{A_1}$ plus lower order terms, the leading
term must be $b(A_1)$ with coefficient 1. If we can prove that the
other coefficients actually are polynomials in $q^2$ without constant
term then the proof follows from Proposition~\ref{4.3}. We proceed to
prove this: First we expand
\begin{equation}
  X_t= Z_{0,0}X_o +\sum_j(-q^2)^{\ell_{\sigma_j}}Z(C_{\sigma_j}),
\end{equation}
where the powers of $q^2$ are strictly positive and where
$Z(C_{\sigma_j})$ is a normalized monomial without contribution from 
$Z_{0,0}$. Thus, $Z(C_{\sigma_j})=b(C_{\sigma_j})+\sum_i\ell_i$,  where
the terms $\ell_i$ are of lower order and have coefficients in $q^2{\mathbb Z}[q^2]$. According to Lemma~\ref{samelem},
\begin{equation}
  (\sum_j(-q^2)^{\ell_{\sigma_j}}Z(C_{\sigma_j}))X_b
\end{equation}
then fulfills the requested condition. It remains to consider $Z_{0,0}X_oX_b$. But here
we notice that $X_oX_b$ again is a member of the dual canonical
basis; $X_oX_b=b(A_3)$ with $A_3$, appropriately cut down, given by \begin{equation}A_3=\left(\begin{array}{ccccc}0&&\cdots&&0\\&2\\\vdots&&\ddots&&\vdots\\&&&2\\0&&\cdots&&1\end{array}\right).\end{equation} This matrix is without contributions from the first row and the
first column. Now notice that
\begin{equation}
  b(A_3)=Z(A_3)+q^2\sum_{G_k<A_3}d_k(q)Z(G_k),
\end{equation}
where the coefficients are in ${\mathbb Z}[q^2]$ and the first row and column in each $G_k$ are zero.  It then follows from
the above remarks that
\begin{equation}
  Z_{0,0} b(A_3)=Z(A_3+E_{0,0})+q^2\sum_{G_k<A_3}d_k(q)Z(G_k+E_{0,0}).
\end{equation}
Expanding the right hand side in terms of the canonical basis, we
get the requested result about the coefficients $c_i$ since the basis change matrix is lower diagonal with
1's in the diagonal and all non-diagonal terms have zero constant
terms. Also notice that $A_3+E_{0,0}=A_1$.

Proposition~\ref{4.3}  implies that $X_tX_b-q^2Y_rY_l$ is bar invariant
and therefore  coincides with the dual canonical basis elements
with the same leading term which is $X_oD$. \qed

\medskip

\begin{Rem} The referee has kindly informed us that another, in a way easier, proof may be obtained as follows: 
Consider the equation
$Z_{1,1}Z_{n,n} = q^2Z_{1,n}Z_{n,1} + \xi^{1,n}_{1,n}$ in ${\mathcal O}_q(M(n))$ and apply the anti-endomorphism $\Gamma$ 
given in \cite[Corollary 5.2.2]{pw} to it. The effect of $\Gamma$ on quantum minors has been computed in
\cite[Lemma 4.1]{kelly}, and from this one can see that the formula is obtained.
\end{Rem}

\bigskip

The following result is a variation of \cite[Lemma~4.1,
Proposition~4.5]{jj}:

\begin{Lem}\label{centerkernel}Let $\Lambda_{L^+}$ be the  $\Lambda$-matrix of
 the Laurent quasi-polynomial algebra generated by the elements of the family ${\mathcal
V}_{L^+}$ and let ${\mathcal C}={\mathcal C}_{L^+}^-$ be the set of $n+m-1$ covariant
quantum minors determined as the variables in ${\mathcal V}_{L^+}$ corresponding
to the points in $L^-$. Let $s=\textrm{corank}(\Lambda_{L^+})$. The kernel of
$\Lambda_{L^+}$ 
is then generated by $s$ elements that have non-zero coefficients at most at the
positions of ${\mathcal C}$.
\end{Lem}

{{}\proof We consider here $m\times n$ matrices; in \cite{jj} it is $n\times r$. We assume $n\leq m$ and $s>0$. Then $s=g.c.d.(m,n)$. Furthermore, if  $m=y\cdot s$, and $n=x\cdot s$, then $x$ and $y$ are odd. Let us label  the covariant elements by integers $j=1,2,\dots ,n+m -1$ corresponding to  the broken line starting at $(1,n)$, passing through $(m,n)$, and terminating at $(m,1)$. The covariant element $\Psi_j$  then has its lower right corner at the $j$th point on the line. Set $\Psi_0=1$. For each $a=1,\dots,s$ the element in the Laurent quasi-polynomial algebra generated by the covariant elements, $$\Phi_a=\prod_{\ell=-x}^{y-1}\Psi_{a+\ell\cdot s+n-1}^{(-1)^\ell},$$ 
commutes with all elements $Z_{i,j}$, and hence with all elements in ${\mathcal V}_{L^+}$. The proof of this is in \cite{jj}. Using Lemma~\ref{2.22} one  easily constructs elements $v(\Phi_1),\cdots,v(\Phi_s)$ in the kernel of $\Lambda_{L^+}$. The element $v(\Phi_a)$ is given with zeros everywhere except at the covariant elements $\Psi_k$ where the coordinate is $\sum_{\ell=-x}^{y-1}(-1)^\ell \delta_{k.a+s\cdot\ell+n-1}$. As these elements evidently are linearly independent, the claim follows. \qed}

\bigskip

\begin{Thm}\label{seed1}Consider two broken lines $L_1< L_2$ in 
$M_{m,n}({\mathbb Z}_+)$ such that $L_2$ is a closest bigger line to $L_1$. Let
${\mathcal Q}_{L_1}=({\mathcal V}_{L_1}, \Lambda_{L_1},B_{L_1})$ and ${\mathcal
Q}_{L_2}=({\mathcal V}_{L_2}, \Lambda_{L_2},B_{L_2})$ be quantum seeds
corresponding to these such that the set of non-mutable elements in both cases
is ${\mathcal C}={\mathcal C}_{L^+}^-$. Then ${\mathcal Q}_{L_1}$ can be
obtained from ${\mathcal Q}_{L_2}$ through a sequence of moves which alternate
between modifications to the $B$ matrix at a given step and quantum mutations.
The modifications to the $B$ matrices do not affect their principal parts.

\end{Thm}

\proof Let $L_2\neq L^-$ be an otherwise  arbitrary broken line and let $(i,j)$
be a concave
corner of $L_2$, specifically, assume that $(i,j+1), (i,j),(i+1,j)$ are 
points on
the broken line $L_2$. If we replace $(i,j)$ by $(i+1, j+1)$ while keeping the
other points, we get a broken line
$L_1<L_2$  such that $L_2$ is a closest bigger line to $L_1$. Consider
quantum seeds ${\mathcal Q}_{L_2}$ and ${\mathcal Q}_{L_1}$ as in the
statement of the theorem  corresponding to this configuration. We
will
construct a sequence of interim quantum seeds ${\mathcal Q}_a=({\mathcal V}_{a},
\Lambda_{a},B_{a})$  and $\tilde{{\mathcal Q}}_a=({\mathcal V}_{a},
\Lambda_{a},\tilde{B}_{a})$ such that
\begin{equation}\label{smut}{\mathcal Q}_{L_2}={\mathcal Q}_{0}\rightarrow
\tilde{{\mathcal Q}}_{0}\Rightarrow{\mathcal
Q}_{1}\rightarrow \tilde{{\mathcal
Q}}_{1}\Rightarrow\dots\rightarrow \tilde{{\mathcal Q}}_{m-i}={{\mathcal
Q}}_{L_1} .\end{equation}

The double arrows are quantum mutations while the single arrows indicate some
change on the level of the B-matrix.

Without loss of generality, we may assume that  $j+m-i\leq n$. By
construction, the quantum minor $\xi^{\{i, i+1,\cdots, m\}}_{
\{j, j+1,\cdots, j+m-i\}}$ is both a quantum cluster variable for the
quantum seeds associated to $L_2$ and  to $L_1$, but labeled by
different points, namely, labeled by $(m, j+m-i)$ in the quantum
seed associated to $L_2$ and labeled by $(i, j)$ in the quantum
seed associated to $L_1$. This quantum minor is not affected by the
following
mani\-pulations. The  quantum minors $\xi^{\{i\}}_{\{j\}},\xi^{\{i,i+1\}}_{
\{j,j+1\}},\cdots, \xi^{\{i, \cdots, m-1\}}_{\{j,\cdots,
j+m-1-i\}}$ are changed, step by step, into $\xi^{\{i+1\}}_{
\{j+1\}},\xi^{\{i+1,i+2\}}_{\{j+1,j+2\}},\cdots, \xi^{\{i+1, \cdots, m\}}_{
\{j+1,\cdots, j+m-i\}}$, and all other quantum cluster variables stay unchanged.

Specifically, we do the following sequence of replacements:
\begin{equation}\label{replace}\begin{array}{ccc}
\xi^{\{i\}}_{\{j\}}&\mapsto& \xi^{\{i+1\}}_{\{j+1\}}\\
\xi^{\{i,  i+1\}}_{\{j, j+1\}}&\mapsto
  & \xi^{\{i+1,  i+2\}}_{\{j+1, j+2\}}\\
&\vdots&\\
\xi^{\{i, \cdots, i+a\}}_{\{j,\cdots, j+a\}}&\mapsto&\xi^{\{i+1, \cdots,
i+1+a\}}_{\{j+1,\cdots, j+1+a\}}\\
&\vdots&\\\xi^{\{i, \cdots, m-1\}}_{\{j,\cdots,
j+m-1-i\}}&\mapsto&\xi^{\{i+1, \cdots, m\}}_{
\{j+1,\cdots, j+m-i\}}
\end{array}\quad .\end{equation}
We claim that each replacement is a quantum
 mutation in the sense of Berenstein, Zelevinsky. Indeed, quantum mutations are
determined by $B$ matrices and we know, modulo the kernel of $\Lambda_a$, {{}the relevant column of $B_a$ by
Corollary~\ref{level}}: At any level $a<m-i$, the quantum minors
 $\xi^{\{i,\cdots, i+a\}}_{\{j,\cdots,j+a\}}=X_t^{(a)}$, $\xi^{\{i+1,\cdots,
i+1+a\}}_{\{j+1,\cdots,j+1+a\}}=X_b^{(a)}$,  $\xi^{\{i,\cdots,
i+1+a\}}_{\{j,\cdots,j+1+a\}}=D^{(a)}$, $\xi^{\{i+1,\cdots,
i+a\}}_{\{j+1,\cdots,j+a\}}=X_o^{(a)}$,
$\xi^{\{i+1,\cdots, i+1+a\}}_{\{j,\cdots,j+a\}}=Y_l^{(a)}$, and
$\xi^{\{i,\cdots,
i+a\}}_{\{j+1,\cdots,j+1+a\}}=Y_r^{(a)}$ constitute an ${\mathcal M}$-set. The
elements $X_t^{(a)}, D^{(a)}$,  $X_o^{(a)}, Y_l^{(a)}$, and  $Y_r^{(a)}$ are all
 quantum cluster variables in ${\mathcal V}_a$.  On the other hand, for
$a\geq1$, $X^{(0)}_t,\dots, X^{(a-1)}_t$ are not variables in ${\mathcal V}_a$.
{{} By easy checking of which variables from ${\mathcal V}_a$ contain
one or the other of $Z_{i,j}$ and $Z_{i+a+1,j+a+1}$, or
both, it follows from Corollary~\ref{level} that the
element}
$X_o^{(a)}D^{(a)}(Y_l^{(a)})^{-1}(Y_r^{(a)})^{-1}$
 commutes with all the quantum cluster variables in ${\mathcal V}_a$ except $X_t^{(a)}$. 
 In fact,
$$X_t^{(a)}X_o^{(a)}D^{(a)}(Y_l^{(a)})^{-1}(Y_r^{(a)})^{-1}=q^{-4}X_o^{(a)}D^{
(a)}(Y_l^{(a)})^{-1}(Y_r^{(a)})^{-1}X_t^{(a)}.$$
 By Lemma~\ref{2.22} this implies that, modulo the kernel of $\Lambda_a$,  the only non-zero entries
in the column of $B_a$ corresponding
to the variable $X_t^{(a)}$ are at the row positions of the
variables $Y_l^{(a)}$ and $Y_r^{(a)}$ where it is $1$, and at the row positions
of
the variables $X_o^{(a)}$ and $D^{(a)}$ where it is $-1$. We then change
(if needed) $B_a$ into $\tilde{B}_a$ such that the column in the latter
corresponding to $X_t^{(a)}$ is non-zero precisely at the mentioned 4 places.
Notice that this change only involves the kernel of $\Lambda_a$.  By considering $B_a^T\Lambda_a$ it follows that any element in the kernel has non-zero coefficients at most at the places of the non-mutable elements ${\mathcal C}$. With this,    {{} the quantum
 mutation of $X_t^{(a)}$ to some new element $\left(X_t^{(a)}\right)^\prime$  in the sense of (\ref{muta})} can be performed;  and by combining Lemma~\ref{red1} with
Theorem~\ref{4.4}, it follows that $X_b^{(a)}=\left(X_t^{(a)}\right)^\prime$ is the target of this mutation. We then perform
this quantum mutation and obtain a new interim quantum seed ${\mathcal
Q}_{a+1}=({\mathcal V}_{a+1}, \Lambda_{a+1},B_{a+1})$. The variables
$Y_l^{(a)}, Y_r^{(a)}, a=0,1,\dots$  are both in ${\mathcal V}_{a}$ and  in
${\mathcal V}_{a+1}$. In the step $a+1$, $D^{(a)}=X_t^{(a+1)}$ and, most
importantly,
$X_o^{(a+1)}=X_b^{(a)}$ which is now a variable in ${\mathcal V}_{a+1}$. In this
way we
can carry out the entire transition from
$L_2$ to $L_1$.   Hence
our changing of the set of variables for a broken line at a concave point is
obtained through a sequence of steps alternating between quantum
mutations in the sense of Berenstein and Zelevinsky, and changes to the $B$
matrix.  

\smallskip

\qed

\medskip

[As an aside, we observe that we, starting at the top,  could break off the
above
replacements at any lower level, but we shall not find it useful to do so.]

\medskip

\begin{Rem}
 In case $\Lambda_{L_2}$ is invertible all modifications to the $B$ matrices in Theorem~\ref{seed1} are
trivial. {{}Indeed,  a modification makes changes involving only the kernel of $\Lambda_{L_2}$.} Hence, in this case the two quantum seeds are equivalent by quantum
mutations. If we can embed  our algebra
${\mathcal O}_q(M(m,n))$ into a bigger algebra ${\mathcal O}_q(M(m_1,n_1))$
such that ${\mathcal Q}_{L_2}$ and ${\mathcal Q}_{L_1}$ are the restrictions of
quantum seeds ${\mathcal Q}^E_{L_2}=({\mathcal V}^E_{L_2},
\Lambda^E_{L_2},B^E_{L_2})$ and ${\mathcal Q}^E_{L_1}=({\mathcal V}^E_{L_1},
\Lambda^E_{L_1},B^E_{L_1})$ with $\Lambda^E_{L_2}$ invertible, then we need at
most make modifications to $B_{L_2}$ and  $B_{L_1}$.

We shall see later that this can always be accomplished. 
\end{Rem}

\medskip

Now, for
each broken line $L$, we have a family of $mn$ $q$-commuting 
quantum minors which by construction is a generating set of the fraction field of
the Noetherian domain ${\mathcal O}_q(M(m,n))$.

\begin{Cor}Let $L$ be an arbitrary broken line and let
$\xi_I^J$ be any solid quantum minor. Then $\xi_I^J$ can be
written as a $q$-Laurent polynomial with coefficient in ${\mathbb
Z}_+[q^{2},q^{-2}]$ of the cluster variables associated to
$L$.\end{Cor}

\proof By our construction using broken lines, one can see that the
solid  quantum minor $\xi_I^J$ belongs to some quantum seed
associated to a broken line $L^\prime$. By the
above theorem, $\xi_I^J$ can be obtained through a sequence of
quantum mutations from the quantum cluster variables associated to
$L$. Now the statement follows
from the quantum Laurent phenomenon established in \cite[Corollary~5.2]{bz}. \qed

\medskip

\begin{Rem} Recall that a real matrix $A$ is totally positive (resp.
 totally non-negative) if all of its minors are positive (resp.
non-negative). In \cite{GP}, it is shown that a matrix is totally
positive if all of its solid minors are positive. Moreover, in
\cite{FZP}, it is shown that a matrix is totally positive if some
specially chosen minors (in fact a cluster) are positive. The above
result is related to the totally positivity of real matrices.
Specializing $q$ to $1$, we obtain a family of seeds (associated to
broken lines) which are mutation equivalent to each other. To test if a
matrix is totally positive one only needs to check if the minors in an
arbitrary cluster associated to a broken line are positive.
\end{Rem}

\medskip

\subsection{Quantum line mutations}

\begin{Def} \label{qmut} In the general setting of ${\mathcal F}_L= {\mathcal
F}_{L^+}$, let $L_1$ be a closest bigger line to the line $L$. Assume
  the configurations are as in Definition~\ref{closest}. The
  restricted quantum line mutation $\mu_R(L_1,L)$ is the map ${\mathcal
V}_{L_1}\rightarrow {\mathcal V}_{L}$ given as the composite map 
  (\ref{replace}) where $(i,j)$ is replaced by $(c-1,d-1)$. {{}Assume that the set of non-mutable elements is ${\mathcal C}_{L^+}^-$. }

  If ${\mathcal Q}_{L_1}=({\mathcal V}_{L_1}, \Lambda_{L_1},B_{L_1})$ and
${\mathcal Q}_{L}=({\mathcal V}_{L}, \Lambda_{L},B_{L})$ are quantum seeds,  the
{\bf quantum line mutation} $\mu(L_1,L):{\mathcal Q}_{L_1}\rightarrow{\mathcal Q}_{L}$ is
a {{}process} as given by the analogue of (\ref{smut}) but where it furthermore is
demanded that at each level $i$, $\tilde{\mathcal Q}_i={\mathcal Q}_i$. The existence of this will be established in Proposition~\ref{bigprop}. For
  practical purposes, we also consider the trivial quantum mutation as a quantum
line
  mutation and denote it by $\mu(L,L)$.
\end{Def}

\medskip

\begin{Def} \label{qmut1} In the general setting of ${\mathcal F}^-_L$, let
$L_2\leq L_3\leq L$ be broken lines such that $L_3$ is a closest bigger line to
the line $L_2$. The
  quantum line mutation $\mu^L(L_3,L_2)$ is the {{}process} ${\mathcal
Q}_{L_3,L} = ({\mathcal V}_{L_3,L}, \Lambda_{L_3,L},B_{L_3,L})\rightarrow {\mathcal
Q}_{L_2,L} = ({\mathcal V}_{L_2,L}, \Lambda_{L_2,L},B_{L_2,L})$ defined in analogy with
Definition~\ref{qmut}. In particular, $\mu(L_1,L)=\mu^{L^+}(L_1,L)$. We denote the inverse of $\mu^L(L_3,L_2)$ by $\mu^L(L_2,L_3)$.
\end{Def}

\medskip

We have the following diamond lemma for quantum line mutations, cf.
(\cite{diamond}):

\begin{Lem}\label{dlemma}Let $L_1\leq L$.
Let $\mu^L(L_1,L_1\downarrow(c_1,d_1))$ and $\mu^L(L_1,L_1\downarrow(c_2,d_2))$
be quantum line mutations. Then $\mu^L(L_1\downarrow(c_1,d_1),
(L_1\downarrow(c_1,d_1))\downarrow (c_2,d_2))$ and
$\mu^L(L_1\downarrow(c_2,d_2), (L_1\downarrow(c_2,d_2))\downarrow (c_1,d_1))$
are quantum line mutations. Furthermore,
\begin{eqnarray*}\mu^L(L_1\downarrow(c_1,d_1),
(L_1\downarrow(c_1,d_1))\downarrow (c_2,d_2))\circ\mu^L(L_1,
L_1\downarrow(c_1,d_1))=\\\mu^L(L_1\downarrow(c_2,d_2),
(L_1\downarrow(c_2,d_2))\downarrow (c_1,d_1))\circ\mu^L(L_1,
L_1\downarrow(c_2,d_2)).
\end{eqnarray*}
\end{Lem}

\smallskip

\proof The key to this Lemma is Corollary~\ref{import} as well as the explicit
formulas (\ref{eq:E-entries}) and (\ref{Zel}). The mutation which  does the
replacement $X_t^{(a)}\mapsto X_b^{(a)}$ makes changes  to the rows in $B$
corresponding to the quantum minors $ D^{(a)}, X_o^{(a)}, Y_l^{(a)}$, and
$Y_r^{(a)}$. It follows by direct inspection that if an entry in the row of
$X_t^{(a)}$ in some position $v$ is zero then the column of $v$ stays unchanged
under the quantum mutation. At this level the entry at the position of
$X_t^{(a)}$ of course is zero. Clearly,  $X_t^{(a)}$ is not a member of any of
the subsequent sets $ D^{(a+p)}, X_o^{(a+p)}, Y_l^{(a+p)}, Y_r^{(a+p)}$, where
$p=1,\dots, p_0$ for some specific positive integer $p_0$. It follows that the
positions in the row of  $X_t^{(a)}$ corresponding to the later values
$X_t^{(a+p)}$ must be zero since we know precisely what the column of
$X_t^{(a+p)}$ looks like. These considerations can easily be extended to include
the case of a second quantum line mutation since none of the variables
$X_t^{(a)}$ take part in any way in the second quantum line mutation. In the
case where $(c_2,d_2)=(c_1-1,d_1+1)$ there is an overlap of variables in the sense
that  the $Y_r$ variables belonging to $(c_1,d_1)$ play the role of $Y_l$
variables belonging to $(c_2,d_2)$, but this is easily taken care of:  They are
not the sources or targets of mutations and then the effects of the two
different quantum line mutations on the rows of such elements are independent of
each other. The crucial observation is that neither of the two quantum line
mutations affects the rows  of the variables involved in the other. \qed

\medskip

The following result concerning independence of paths, follows easily since one may fill in diamonds as in Lemma~\ref{dlemma}:

\begin{Cor}\label{iop}
If $L_1\leq L_2\leq\dots\leq L_{n-1}\leq L_n\leq L$ and $L_1\leq
L'_2\leq\dots\leq L'_{n-1}\leq L_n$ are broken lines such that at each level
the bigger line is a closest bigger line to the neighboring smaller one. Then
$$\mu^L(L_2,L_1)\circ\dots\circ \mu^L(L_n,L_{n-1})=
\mu^L(L'_2,L_1)\circ\dots\circ \mu^L(L_n,L'_{n-1}).$$
\end{Cor}

\medskip

In view of Corollary~\ref{iop} we extend our definition of a quantum line
mutation to the following

\medskip

\begin{Def} Let $L_1\leq\L_n\leq L$ be broken lines. The quantum line mutation
$\mu^L(L_n,L_1)$ is the composite of any sequence as in Corollary~\ref{iop}
between $L_1$ and $L_n$. 

Let $L_a\leq L$ and $L_b\leq L$ be broken lines. The quantum line mutation
$\mu^L(L_a,L_b)$ is defined in terms of  any broken line $L_c\leq L_a,L_b$ as
$$\mu^L(L_a,L_b)=\mu^L(L_b,L_c)^{-1}\circ\mu^L(L_a,L_c).$$
\end{Def}

\bigskip

We shall also need

\begin{Def}A position $(c,d)\in T_L\cup L$ is called attractive with respect to
$T_L\cup L$ if either there exist $i>0,\  j>0$ such that $(c-i,d)\in T_L\cup L$
and $(c,d+j)\in T_L\cup L$ or if there exist $i>0,\  j>0$ such that $(c+i,d)\in
T_L\cup L$ and $(c,d-j)\in T_L\cup L$. Clearly, if $(c,d)$ satisfies the first
condition of attraction then $(c-i,d+j)$ satisfies the second, and vice versa. 
If $(c,d)$ is not attractive we call it repulsive.
\end{Def}

The following is obvious

\begin{Lem} The concave corners of $L$ are repulsive. The point $(m,n)$ is also
repulsive.
\end{Lem}

\bigskip

\subsection{Covariant elements}

We extend Definition~\ref{covdef} in the obvious way to ${\mathcal
    O}_q(T_L\cup L)$. The next observation we wish to make is that the seeds we
construct are minimal in the following sense:

\begin{Prop}\label{cov-prop}
The set of covariant elements for ${\mathcal
    O}_q(T_L\cup L)$ is contained in the sub-algebra generated by the $m+n-1$
elements in ${\mathcal
    C}^-_{L}$.
\end{Prop}

\proof First of all it is clear that the elements in  ${\mathcal
    C}^-_{L}$ are covariant, and hence, so is any monomial in these. 
    
    Since there is a unique smallest element in the set of broken lines, this may be proved
by induction. For the line $L^-$ it is clear that we have a
quasi-polynomial algebra so here, the claim is trivial.
Consider then a line $L$ for which the claim is true and let $L_1$
be a closest bigger line. Assume the configurations are as in
Definition~\ref{closest}. (Thus, $(i,j)=(c-1,d-1)$.) It is clear that ${\mathcal
O}_q(T_L\cup
L_1)$ is obtained by adjoining $Z_{c-1,d-1}$ to ${\mathcal
O}_q(T_L\cup L)$. There is a unique element  $X_b$
from ${\mathcal C}_L^-$ having its upper left corner in $(c,d)$. Depending on the configuration,  $X_b=X_b^{(m-c)}$, or $X_b=X_b^{(n-d)}$. 
This is the largest solid quantum minor with its upper left corner in
this position and completely contained in ${\mathcal O}_q(T_L\cup L)$.
It is clear from Proposition~\ref{q-cor} that this is the only element
from ${\mathcal C}_L^-$ which does not $q$-commute with
$Z_{c-1,d-1}$. On the other hand, when $(c-1,d-1)$ is viewed as an
element in $S_L$, the  variable
$D=D^{(m-c)}$, respectively  $D=D^{(n-d)}$, does, by Proposition~\ref{a-l},
$q$-commute with all the $Z_{i,j}$ in ${\mathcal O}_q(T_L\cup L)$ - and clearly also with
$Z_{c-1,d-1}$. Next observe that evidently $D
\in{\mathcal C}_{L_1}^-$. 

Suppose then that $C\in {\mathcal O}_q(T_{L_1}\cup L_1)$ is covariant. It is
clear that 
$${\mathcal O}_q(T_{L_1}\cup L_1)\subseteq {\mathcal O}_q(T_{L}\cup
L)[D, X_b^{-1}].$$  Both adjoined elements are covariant as far as
 ${\mathcal O}_q(T_{L}\cup L)$ is concerned, and it follows easily that $C$ must
be a polynomial in the variables from ${\mathcal
    C}^-_{L}$ together with $D$ and  $X_b^{-1}$. The element
$Z_{c-1,d-1}$ $q$-commutes with all these generators except $X_b$ and this
easily implies that $X_b^{\pm1}$ cannot appear.
    
    \qed

\medskip

\begin{Prop}Consider the quadratic algebra ${\mathcal
    O}_q(T_L\cup L)\subseteq{\mathcal O}_q(M(m,n))$. Then there is a non-trivial
center if and only if $m=n$ and $L=L^+$. This center is generated by
$\det_q(n)$.
\end{Prop}

\proof Consider the covariant element $M=M_{m,n}\in{\mathcal C}_L^-$ corresponding to the position
$(m,n)$. If $N$ is any other covariant element and $MN=q^\delta NM$ then
$\delta\leq0$. This follows by easy inspection. Hence, since
any central element must be a polynomial in the elements in ${\mathcal C}_L^-$,
the central element must a polynomial  in those elements from ${\mathcal C}_L^-$
that are quantum minors of $M$. If there are elements $Z_{i,j}$  in the algebra not occurring in $M$,
then there will be non-trivial commutation relations between these and the
quantum minors from $M$. Thus, there can be no positions outside $M$ if $M$ is to be central. The remaining
details
now follow from the classical result of Parshall and Wang (\cite{pw}). \qed

\medskip

\bigskip

\section{On compatible pairs}\label{5}

We now settle the existential questions implicitly raised in
Theorem~\ref{seed1}, Definition~\ref{qmut}, and Definition~\ref{qmut1}.

\begin{Prop}\label{bigprop} To a given broken line $L$ in $M_{m,n}({\mathbb Z}_+)$ one can construct the
following data
${\mathcal D}_L$:
\begin{itemize}
\item An ordering of the set of variables ${\mathcal V}_L={\mathcal
    V}^+_L\cup {\mathcal V}^-_L$ of the broken line such that the
  variables in ${\mathcal V}^-_L$ are assigned the numbers from $1$ to
  $N_L$ and the remaining variables the numbers from $N_L+1$ to
  $mn$. Let $c_L^-=\#{\mathcal C}^-_L$. (${\mathcal C}^-_L$ is the set of   non-mutable variables when
  ${\mathcal V}^-_L$ is considered in the ambient space ${\mathcal F}_L^-$. We have, of
  course, that $c_L^-=m+n-1$ but it is convenient to have the additional
  notation.)  We will even assume that the mutable variables are
  assigned the numbers from $1$ to $\tilde N_L=N_L-c_L^-$.
\item The $mn\times mn$ matrix $\Lambda_L$ corresponding to the variables ${\mathcal
    V}_L$ of the broken line. 
\item A matrix $B_L$ of size $mn\times (m-1)(n-1)$
  such that $(\Lambda_L,B_L)$ is a
  compatible pair.
\item The $N_L\times N_L$ matrix $\Lambda^0_L$ of the variables in
  ${\mathcal V}^-_L$. We view this as a sub-matrix of
  $\Lambda_L$. 
  \item  An $N_L\times \tilde N_L$ matrix $B^0_L$ such that
$(\Lambda^0_L,B_L^0)$ is a compatible pair for the set of
variables ${\mathcal V}_L^-$.
If one defines a $mn\times \tilde N_L$ matrix $B^{\textrm{R}}_L$ by adding
$mn-N_L$ rows of zeros to $B_L^0$ such that $B_L^0$ occupies the top rows of
$B_L$, the following holds in addition:
\begin{itemize}
  \item[$\rightarrow$] $\Lambda_LB^{\textrm{R}}_L=-2D_L$, where
$D_L=I_{\tilde
      N_L}\oplus 0$ is the $mn\times \tilde N_L$ matrix consisting
    of $I_{\tilde N_L}$ in the top $\tilde N_L\times {\tilde N_L}$
    corner augmented by an appropriate number of rows of zeros. Here,
    $I_{\tilde N_L}$ is the ${\tilde N_L}\times {\tilde N_L}$ identity
    matrix.
\item[$\rightarrow$] $B^{\textrm{R}}_L$ is a sub-matrix of $B_L$. 
\end{itemize}
\end{itemize}
In the case of  ambient space  ${\mathcal F}^-_{L}$, consider
the quantum seed ${\mathcal Q}^-_L=({\mathcal V}^-_L,\Lambda_L^0,B_L^0)$ with the
set non-mutable elements given as ${\mathcal C}_L^-$. If $L_1$ is a broken line
in  $M_{m,n}({\mathbb Z}_+)$ and $L_1<L$, there exists a quantum seed ${\mathcal
Q}_{L_1,L} = ({\mathcal V}_{L_1,L}, \Lambda_{L_1,L},B_{L_1,L})$ which is
equivalent by the quantum line mutations $\mu^L(L,L_1)$ to ${\mathcal Q}^-_L$ and where the
pairs $(\Lambda_{L_1,L},B_{L_1,L})$ and $(\Lambda^0_{L_1},B^0_{L_1})$ are related in a way that generalizes in an obvious manner the way $(\Lambda_{L_1},B_{L_1})$ is related to $(\Lambda^0_{L_1},B^0_{L_1})$. 
\end{Prop}

\label{b}

\pof A short proof would be to say that this follows by bootstrapping. We give
here a more detailed proof using the same principle: First assume that
$n=m+1$ (or, analogously, $n=m-1$). The existence of $B_L$ will follow from the
first
parts of the
proof.  We prove the claims involving $B^{\textrm{R}}_L$ and $B_L^0$ by
induction on the partial order on the set of broken lines. The
induction starts with the line $L^-$. There are no mutable elements in
${\mathcal V}^-_{L^-}$ so $B^{\textrm{R}}_{L^-}$ and $B^0_{L^-}$ are empty. The
other structure we start with is ${\mathcal V}_{L^-}$ and a compatible pair 
$(\Lambda_{L^-},B_{L^-})$ connected with this set of variables. The set of  non-mutable
elements is ${\mathcal C}={\mathcal C}^-_{L^+}$ as before. The matrix
$\Lambda_{L^-}$ can, naturally, be explicitly written down. It is known from
\cite{jj}[Proposition~4.5] that   $\Lambda_{L^-}$ is
invertible as a real matrix. Furthermore it follows from Proposition~4.11 {{}therein} that
there is a block decomposition into $2\times2$ skew integer matrices. It
follows then  from the discussion on p. 85 in \cite{jj} that there exists an
integer matrix $A$ with $\det A=1$ such
that $A(\Lambda_{L^-})A^t=\tilde D$ where $\tilde D$ is a block diagonal matrix
consisting of $\frac12 m(m-1)$ $2\times 2$ blocks  
$\left(\begin{array}{cc}0&2\\-2&0\end{array}\right)$ and  $m$ $2\times 2$ blocks
$\left(\begin{array}{cc}0&1\\-1&0\end{array}\right)$.  The results in \cite{jj}
are obtained where $q$ is a root of unity. The end results on the structure of $\Lambda_{L^-}$ are, however, independent of special choices of $q$, and hence apply to the block
diagonalization of the integer matrix $\Lambda_{L^-}$ as in our 
case. It follows that $2\Lambda_{L^-}^{-1}=2{A^T}(A\Lambda_{L^-}
A^T)^{-1}A$ is an integer matrix. The existence of
$B_{L^-}$ follows easily from this as the $(m(m+1)\times (m(m+1)-c^-_L)$
sub-matrix of $-2\Lambda_{L^-}^{-1}$ consisting of the first $(m(m+1)-c^-_L)$
columns. At the moment, it is only the existence and uniqueness of $B_{L^-}$ (up
to multiplication by a positive integer) that matters. After we have presented
the induction
step, we encourage the reader to take it right at the start as a simple
exercise.

Suppose now that we have a line $L$ with
data ${\mathcal D}_L$. Let $L_1>L$ be a broken line closest to $L$. We must now
construct the data ${\mathcal D}_{L_1}$ for $L_1$. For this purpose we consider
the inverse, $\mu(L,L_1)$, of
the quantum line mutation $\mu(L_1,L)$. We view these mutations as taking
place inside the full $m(m+1)\times m(m+1)$ matrix algebra, but we shall
keep a keen eye on its relation to the ambient spaces ${\mathcal
F}^-_L$ and ${\mathcal F}^-_{L_1}$. Specifically, how ${\mathcal
Q}^-_{L_1}$ grows from ${\mathcal Q}^-_L$.

The mutation $\mu(L,L_1)$ begins with a mutation of the form of the
reverse of the bottom line in (\ref{replace}) and where the element we
mutate from, $X_b^{(m-i-1)}$,  is
a covariant, viz.  non-mutable, element of ${\mathcal V}^-_L$. This is not
represented in the matrices $B^0_L$ and $B_L^R$, so we define a new matrix
$\tilde B^R_{L}$ by joining one new column  $c_{b,m-i-1}=c(X_b^{(m-i-1)})$ to
$B^R_L$ in the position $\tilde
N_L+1$ and labeled by $X_b^{(m-i-1)}$. Denote by $\tilde B^0_L$ the sub-matrix
of $\tilde B^R_L$ with the same column numbers and having rows corresponding to
those of $B^0_L$ together with an additional row labeled by $D^{(m-i-1)}$. If
the approach is possible, it follows
from Corollary~\ref{import} in combination with
Definition~\ref{def:compatible-triple} and Lemma~\ref{2.22}  what the
added column
must look like: There should be the value -1 at the
positions corresponding to $X_o^{(m-i-1)}$ and $D^{(m-i-1)}$ and the value 1 at
the positions of
$Y_r^{(m-i-1)}$, and $Y_l^{(m-i-1)}$. All other entries must be zero.
 $\Lambda_{L}(c_{b,m-i-1})$ is a column with exactly one non-zero
coefficient $-4$. This occurs at the position of $X_b^{(m-i-1)}$.
We have that $X_o^{(m-i-1)},
Y_r^{(m-i-1)}$, and $Y_l^{(m-i-1)}$ are variables of ${\mathcal  V}^-_L$. The
element $D^{(m-i-1)}$ is a variable in ${\mathcal V}_L$ but is not a variable of
${\mathcal V}^-_L$. It is a covariant (non-mutable) element of ${\mathcal
V}^-_{L_1}$.
To begin with we consider the set of variables  ${\mathcal V}^-_{L,D}={\mathcal V}^-_{L}\cup
\{D^{(m-i-1)}\}$. Corresponding to this
we have a $\Lambda$-matrix $\tilde\Lambda^0_{L}$ with one more
column and row than $\Lambda^0_L$. Since $-\frac12B_L$ is a part of the
inverse matrix of $\Lambda_L$ it is clear that the
matrix $\tilde B^R_{L}$ is part of the matrix $B_{L}$. Furthermore, by
construction, $(\tilde \Lambda^0_L,\tilde B^0_L)$ is a compatible pair for  ${\mathcal V}^-_{L,D}$. It is also obvious, using the equations involved in (\ref{replace})  in the proof of  Theorem~\ref{seed1}, that the skew field of fractions generated by ${\mathcal V}^-_{L,D}$ is exactly ${\mathcal F}^-_{\Lambda_1}$. We now
perform the mutation inside the full matrix algebra at the variable $X_b^{(m-i-1)}$. On the level of variables, with the given column of $B_{L}$, this is, according to Theorem~\ref{4.4}, exactly the mutation

$$X_b^{(m-i-1)} \mapsto X_t^{(m-i-1)}.$$

Correspondingly, we obtain an interim  pair
$(\Lambda_{L}^{(m-i-1)},B_{L}^{(m-i-1)})$ which corresponds  to an interim set of variables ${\mathcal V}_{L,(m-i-1)}$.  At the same time we perform a mutation in ${\mathcal F}^-_{L_1}$ at the same variable but using the pair $(\tilde \Lambda^0_L,\tilde B^0_L)$. Here we obtain an interim  set of variables ${\mathcal V}^-_{L,(m-i-1)}$ which is a part of ${\mathcal V}_{L,(m-i-1)}$. The two mutations do not differ in what they do to the variable $X_b^{(m-i-1)}$.  The difference lies entirely on the level of compatible pairs. The parts of the mutations which involve the matrix $E_i$ in (\ref{eq:E-entries}) affect the rows corresponding to the variables $X_o,D$ and we get the same contributions for both mutations as far as  $\tilde B^R_{L}$ is concerned. It is less transparent at first what happens at the place of the matrix  
$F_i$. Here there are two different matrices for the two mutations. The  
difference relies on the fact that $D^{(m-i-1)}$, and one or both of
$Y_r^{(m-i-1)},Y_l^{(m-i-1)}$ are
mutable in the full algebra but non-mutable in the
small algebra. At this moment of the proof we are only concerned with what happens to the
columns of $\tilde B^R_{L}$ inside $B_{L}^{(m-i-1)}$. The remaining columns of $B_{L}^{(m-i-1)}$, though,  are known in
principle. The two mutations differ only in  what happens to the columns outside $B_{L}^{(m-i-1)}$ where, according to (\ref{Zel}), a  multiple of the column
corresponding to $X_b$  may be added. Furthermore, if one or both of $Y_r^{(m-i-1)},Y_l^{(m-i-1)}$ are
non-mutable in
${\mathcal V}^-_{L}$ they stay so in ${\mathcal V}^-_{L,D}$. So, since the columns of the non-mutable variables of ${\mathcal V}^-_{L,D}$ are not under discussion, it is clear that the analogous new compatible pair
$(\Lambda^0_{L},B^0_{L})^{(m-i-1)}$, obtained from $(\tilde \Lambda^0_L,\tilde B^0_L)$ by mutation in the
ambient space ${\mathcal F}^-_{L_1}$, is a sub-pair of
$(\Lambda_{L_1},B_{L_1})$.  Now perform the remaining mutations in the
quantum line mutation. These only involve mutable variables and are easily
seen to preserve the general form. Finally, one can reshuffle the
variables to obtain the wanted ordering. Again, this does not change
the general form. Thereby the induction step is completed.

In this way we build up  $B$ matrices with more and
more columns. In the end we reach ${\mathcal Q}_{L^+}=({\mathcal V}_{L^+},\Lambda_{L^+},B_{L^+})$ of the
extremal line $L^+$. Once we have that, we can mutate back, by quantum line mutations, to any quantum seed ${\mathcal Q}_{L}=({\mathcal V}_{L},\Lambda_{L},B_{L})$. We can also stop the growing process at an earlier point, where we have obtained a seed ${\mathcal Q}^-_{L}=({\mathcal V}^-_{L},\Lambda^0_{L},B^0_{L})$ and use mutations $\mu^L$ inside the ambient space ${\mathcal F}^-_L$ to obtain the seeds ${\mathcal
Q}_{L_1,L}$ mentioned in the proposition. Finally recall the independence of path result Corollary~\ref{iop}.

\smallskip

The proof applies to any pair $(m,n)$ for which $\Lambda_{L^-}$ is invertible. Let us then consider the general situation of ${\mathcal O}_q(M(m,n))$. Suppose
for simplicity that $m=n+r$ with $r\geq 2$. We can view this as the sub-algebra 
of ${\mathcal O}_q(M(m,n+r+1))$ generated by the elements $Z_{i,j}$ with $1\leq
i\leq m$ and $2+r\leq j\leq n+r+1$. Any broken line
$L:(1,n)\rightarrow\dots\rightarrow(m,1)$ in $M_{m,n}({\mathbb Z}_+)$ is
similarly considered as a line $\tilde L
:(1,n+r+1)\rightarrow\dots\rightarrow(m,r+2)$ in $M_{m,m+1}({\mathbb Z}_+)$
which is then extended by the segment $\tilde L\rightarrow (m,1)$. This
corresponds to adding the non-mutable covariant variables
$Z_{m,1},\dots,Z_{m,r+1}$ to all sets of variables in all quantum seeds. If we
stipulate that the mutations and other operations in  ${\mathcal
O}_q(M(m,n+r+1))$ should never involve these we clearly get the result as a
sub-case of the full case based on $(m,m+1)$. Finally, the case $m=n$ follows by analogous considerations.

\qed

\medskip

\begin{Rem}
Also for the remaining mutations in the quantum line mutation $\mu(L,L_1)$ we
can write down explicitly the values in the $B$-column which we mutate
from simply by using Corollary~\ref{import} repeatedly. In this way
one can in fact ``explicitly'' write down the compatible pairs at each
step.
\end{Rem}

\medskip

\section{The quantum (upper) cluster algebra  of a  broken line}

We define now some algebras connected with a broken line. Our terminology may
seem a bit unfortunate since the notions of a cluster algebra and an
upper cluster algebra already have been introduced by Berenstein and Zelevinsky in
terms of all mutations. We only use quantum line mutations which form a proper
subset of the set of all quantum mutations. However, it will be a corollary to what
follows that the two notions in fact coincide, and for this reason we
do not introduce some auxiliary notation.

\begin{Def}
The  cluster algebra ${\mathcal A}^-_L$ connected with a broken
line $L$ in $M_{m,n}({\mathbb Z}_+)$ is the ${\mathbb Z}[q]$-algebra generated in the space ${\mathcal F}^-_L$ by the inverses of the non-mutable elements ${\mathcal C}_L^-$ together with the union of the sets of  all  variables
obtainable from the initial seed ${\mathcal Q}^-_L$ by composites of quantum
line mutations $\mu^L(L,L_1)$ with $L_1\leq L$. 
\end{Def}

\medskip
Observe that we include ${\mathcal C}_L^-$ in the set of variables.
\medskip

\begin{Def}
The upper cluster algebra ${\mathcal U}^-_L$ connected with a broken
line $L$ in $M_{m,n}({\mathbb Z}_+)$ is the ${\mathbb Z}[q]$-algebra in  ${\mathcal F}^-_L$ given as the intersection of  all the Laurent quasi polynomial algebras of the sets of variables
obtainable from the initial seed ${\mathcal Q}^-_L$ by composites of quantum
line mutations $\mu^L(L,L_1)$ with $L_1\leq L$. 
\end{Def}

\medskip

\begin{Rem}The results we obtain below are independent of the pairs $(\Lambda,B)$  entering into the quantum seeds.  What enters into the proofs are mutations as in Theorem~\ref{seed1}.
\end{Rem}

\medskip

\begin{Rem}
The algebras ${\mathcal A}^-_L$ and ${\mathcal U}^-_L$ of a broken line $L$ are defined in terms of
some ${\mathcal O}_q(M(m,n))$, but of course, if the line has a segment
$(m,1)\leftarrow (m, u)$, with $(m,u)$ denoting a corner, and $u>1$, then
the elements $Z_{m,1},\dots, Z_{m,u-1}$ are all covariant. Thus,  ${\mathcal O}_q(T_L\cup L)={\mathcal O}_q(T_{L_1}\cup L_1)[Z_{m,1},\dots, Z_{m,u-1}]$ and  ${\mathcal
A}^-_L={\mathcal A}^-_{L_1}[Z_{m,1}^{\pm1},\dots, Z_{m,u-1}^{\pm1}]$, where $L_1$ is what remains
of $L$ after these elements have been removed. Similarly with segments
$(1,n)\rightarrow(u,n)$. In the same spirit, covariant elements may be added if
it is
convenient to view ${\mathcal A}^-_L$ as a part of a quantum cluster algebra based
on some other  ${\mathcal O}_q(M(m_1,n_1))$ with $m\leq m_1$ and $n\leq n_1$. See also the last part of the proof of Proposition~\ref{bigprop}. 
\end{Rem}

It is clear that ${\mathcal O}_q(T_L\cup L)\subseteq {\mathcal U}^-_L$
and that $Y_i^{\pm1}\in {\mathcal U}^-_L$ for all $Y_i\in{\mathcal
C}_L^-$. Indeed, by the $q$-Laurent Phenomenon (\cite[Corollary~5.2]{bz}),   ${\mathcal
O}_q(T_L\cup L)\subseteq {\mathcal A}^-_L\subseteq {\mathcal U}^-_L$.

\medskip

\begin{Thm}\label{mainthm}Let ${\mathcal C}_L^-=\{Y_1,\dots,Y_s\}$. Then,
$${\mathcal U}^-_L={\mathcal O}_q(T_L\cup
L)[Y_1^{\pm1},\dots,Y_s^{\pm1}]={\mathcal A}^-_L.$$
\end{Thm}
We need only establish that ${\mathcal U}^-_L\subseteq{\mathcal O}_q(T_L\cup
L)[Y_1^{\pm1},\dots,Y_s^{\pm1}]$. We will  in the proof of that use the
following

\begin{Prop}\label{prime}
A quantum minor $M\in{\mathcal C}_L^-$ generates a completely prime ideal of ${\mathcal O}_q(T_L\cup L)$. Specifically, it  satisfies the following
crucial property:
 $$\textrm{ If }p_1M=p_2p_3 \textrm{ in } {\mathcal O}_q(T_L\cup L) \textrm{
then
}$$

 $$p_2= p_4M\textrm{ or } p_3= p_5M \textrm{ for some }p_4,p_5 \textrm{ in }
{\mathcal O}_q(T_L\cup L) .$$
\end{Prop}  

\proof  Goodearl and Lenegan proved in  \cite[Theorem~2.5]{gl} that the determinantal
ideal is prime. We can reduce our case, in which $M$ is a covariant
quantum minor,  to
theirs by using a PBW basis of the full set of variables {{}$\{Z_{ij}\}$} in which the variables
of the rows and columns of $M$, henceforth referred to as the variables of $M$, are written to the right. The elements $p_2$ and
$p_3$ may then be written as sums of polynomials in the variables of $M$  with
coefficients (to the left) that are monomials in the variables that are not variables of  $M$. Let
us be specific and say that $M=\xi^{\{i,i+1,\dots,m\}}_{\{j,j+1,\dots,j+m-i\}}$.
Let us order the monomials in the variables not in $M$ so that the points with column number less than $j$
are biggest, and ordered lexicographically with the biggest being the point with
 smallest row and column number. The finer details are irrelevant. Next in the
ordering we take those points having a column number between $j$ and $j+m-i$
with a similar lexicographical ordering. Finally we take those with a column
number bigger than $j+m-i$. Here we chose an opposite ordering. We can then
focus on the monomials that are the biggest in $p_2$ and $p_3$. The point of
the chosen ordering is that one does not pick up bigger terms via (\ref{344})
while rewriting a product. Let  $v_2^0p_2^0$
be the summand in $p_2$ corresponding to the biggest monomial $v_2^0$. Here
$p_2^0$ is a polynomial in the variables of $M$. Let $v_3^0p_3^0$ be the
analogous summand for $p_3$. Regrouping in the product $p_2p_3$ according to our
total ordering results in a unique highest term (up to a factor of $q$ to some
power) $wp_2^0p_3^0$, where $w$ is the highest order element of $v_2^0v_3^0$.
This must then match a term  $v_1^0p_1^0M$ in $p_1M$. Specifically, $v_1^0=w$. By (\cite{gl})],  $p_2^0=p_4^0M$ or $p_3^0=p_5^0M$. Say it is
$p_2^0=p_4^0M$. Since $M$ is covariant with respect to all the variables of $
{\mathcal O}_q(T_L\cup L)$ we can just  drop the expression $v_2^0p_2^0$ from
$p_2$. Indeed, by looking at the biggest elements, we can assume from the
beginning that neither $p_2$ nor $p_3$
contains a summand of the form $pM$ and then argue by contradiction.  \qed

\proof [Proof of Theorem~\ref{mainthm}] We prove this by induction.
For the unique smallest line $L^-$ the algebra ${\mathcal O}_q(T_L\cup
L^-)={\mathcal O}_q (L^-)$ is generated by the covariant elements in
${\mathcal C}_{L^-}^-$. The algebra is quasi-polynomial and there are
no quantum line mutations except the trivial one. Thus the claim is trivially
true. [Actually, there is also a unique line $L_1$ closest to $L^-$
and the situation here essentially corresponds to ${\mathcal O}_q(M(2,2))$. This
case is also true and well-known.] Let us then consider a line $L$
and let $L_1$ be a closest line with $L<L_1$. Let the notation be as
in the proof of Proposition~\ref{cov-prop}. Then  ${\mathcal
C}_{L_1}^-$ is obtained from ${\mathcal C}_{L}^-$ by replacing
$X_b$ by $D$. Suppose that
$u\in {\mathcal U}^-_{L_1}$. Since ${\mathcal U}^-_{L_1}$ is an algebra
it is clear that we may assume that when $u$ is expressed as a
$q$-Laurent polynomial of some set of variables, all powers of the
covariant elements {{} ${\mathcal C}^-_{L_1}$ are non-negative. Then it is true for all allowed sets of variables as in Definition~8.2.}  Moreover,  $L$ is obtained
from $L_1$ by a quantum line mutation and all subsequent quantum line mutations of $L$ are thus
also quantum line mutations of $L_1$. In all these mutations $D$ stays
non-mutable. For all lines $L_2\leq L$, the algebra generated by the  variables from ${\mathcal V}_{L_2,L_1}$ and their inverses is contained in the algebra generated by the  variables from ${\mathcal V}_{L_2,L}$ and their inverses together with $D^{\pm1}$. It is then clear that  ${\mathcal
  U}^-_{L_1}\subseteq {\mathcal U}^-_L[D^{\pm1}]$. Now, the non-mutable
variables  of $L$ are the same as those of $L_1$ with the exception of $X_b$. By
the argument about the positivity of the non-mutable variables we can
then assume
that all these variables except the latter occur with a non-negative power when $u$ is expanded in one of the allowed quasi Laurent algebras. Thus, we can
assume $u\in{\mathcal U}^-_L[X_b^{\pm1}, D]$. By the induction
hypothesis we then have
  \begin{equation}
	u\in{\mathcal O}_q(T_L\cup L)[X_b^{-1}, D],
	\label{u}
\end{equation}
and to meet our goal, we only need to be concerned about the elements with a
strictly negative power of  $X_b$ in each summand.

Naturally, ${\mathcal O}_q(T_L\cup L)$ can be viewed as a sub-algebra of
${\mathcal
O}_q(T_{L_1}\cup L_1)$.
 
Let us denote the initial variables of $L_1$  by $D=X_t^{(m-c+1)},
X_t^{(m-c)},\dots,X_t^{(0)}
$,  $W_1,\dots,W_N$. The initial variables of $L$ are then 
$X_b=X_b^{(m-c)}, X_b^{(m-c-1)}$, $\dots$, $X_b^{(0)}$ $\dots$ $,
W_1,\dots,W_N$. Let us look
at the element $u$. 
This can be written as a $q$-Laurent polynomial in the given initial variables of
$L_1$, one of which is $D$: 

$$u=\sum_{\underline{\alpha},\underline{\beta}} c_{\underline{\alpha},\underline{\beta}}W^{\underline{\beta}}
\prod_{i=0}^{m+1-c}\left(X_t^{(m-c+1-i)}\right)^{\alpha_i}.$$

We can factor out the biggest non-positive powers such that 
\begin{equation}\label{first}u=p_{top}\cdot W^{-\underline{\beta}^0}
\prod_{i=0}^{m+1-c}\left(X_t^{(m-c+1-i)}\right)^{-\alpha^0_i},\end{equation}
where $p_{top}\in {\mathcal O}_q(T_{L_1}\cup L_1)$ and, in particular, $p_{top}$
contains no overall  factor of $D$. We wish to argue by
contradiction and thus assume that the multi-indices
${\underline{\alpha}^0},{\underline{\beta}^0}$ are non-negative, and at least
one $\alpha^0_r$ or $\beta^0_s$ is positive.

Set $Z=Z_{c-1,d-1}$ Then $D=ZX_b$  modulo ${\mathcal O}_q(T_L\cup L)$. By
(\ref{u}) we have

\begin{equation} \label{neg}
u=(\sum_i Z^ip_iX_b^{k_i})X_b^{-\rho},
\end{equation}
where $\forall i: 0\leq k_i<\rho$ and where, furthermore, the elements $p_i\in {\mathcal
O}_q(T_L\cup L)$ are neither divisible by  $Z$ nor by $X_b$. Combining (\ref{first}) and (\ref{neg}), we have

$$(\sum_i Z^ip_iX_b^{k_i})X_b^{-\rho}=p_{top}\cdot
W^{-\underline{\beta_0}}\prod_{i=0}^{m+1-c}\left(X_t^{(m-c+1-i)}\right)^{
-\alpha^0_i}.$$

Now,  in the $q$-Laurent algebra we clearly get

\begin{equation}\label{z-eq}(\sum_i
Z^ip_iX_b^{k_i})\prod_{i=0}^{m+1-c}\left(X_t^{(m-c+1-i)}\right)^{\alpha^0_i}W^{
\underline{\beta}^0}=q^{2\gamma}p_{
top}X_b^\rho,
\end{equation}
where $q^{2\gamma}$ is an irrelevant factor stemming from the $q$-commutativity between $X_b$ and the elements $W_i$. We ignore this
and similar factors in the following.

The first crucial observation is that by Proposition~\ref{prime}, $\alpha_0^0=0$ since
the right hand side of (\ref{z-eq}) clearly does not contain a positive power of $D$.

The next important fact is that the position $(c-1,d-1)$ is repulsive with respect to $T_{L_1}\cup L_1$. This
implies that it is straightforward to look at the highest order terms of $Z$ in
(\ref{z-eq}). In the right hand side we simply write $p_{top}=Z^Su_S+\ell$ where $\ell$ is of lower order, and
where $u_S$ is a polynomial in the variables of $T_L\cup L$. In the left hand
side, let us say that $Z^Kp_KX_b^{k_K}$ is the term containing the highest $Z$
exponent
$K$. We then get additional $Z$ terms from
$\prod_{i=1}^{m+1-c}\left(X_t^{(m-c+1-i)}\right)^{\alpha^0_i}$, and
here the highest $Z$ term is
$Z^{\alpha_S^0}\prod_{i=1}^{m-c}\left(X_b^{(m-c-i)}\right)^{\alpha^0_i}$ where 
$\alpha^0_S=\sum_{i=1}^{m-c+1}\alpha^0_i$. 

Using the repulsiveness again, we get

$$Z^{K+\alpha^0_S}p_K\prod_{i=1}^{m-c}\left(X_b^{(m-c-i)}\right)^{\alpha^0_i}
X_b^{k_K}
=Z^Su_SX_b^\rho.$$

And thus, 

$$p_K\prod_{i=1}^{m-c}\left(X_b^{(m-c-i)}\right)^{\alpha^0_i}=u_SX_b^{\rho-k_K}
.$$

This identity holds in ${\mathcal O}_q(T_L\cup L)$. Since $\rho-k_k>0$ it
follows by Proposition~\ref{prime} 
that $X_b$ must be a right divisor of one of the terms on the left hand side.
The
$X_b^{(a)}$ with $a=0,1,\dots, m-c-1)$ terms of course are impossible in this
respect. Thus, $p_K=\hat{p}X_b$ for
some $\hat{p}$. This is a contradiction to the way  $p_K$ was defined. Hence,
there can be no negative power $X_b^{-\rho}$ in (\ref{neg}). Thus, $u\in{\mathcal O}_1(T_L\cup L)[D]\subseteq {\mathcal O}_q(T_{L_1}\cup L_1)$.

\qed

\medskip

Since we have established the more restrictive inclusion ${\mathcal U}^-_L\subseteq{\mathcal A}^-_L$ we get 

\begin{Cor}The algebras ${\mathcal A}^-_L$ and ${\mathcal U}^-_L$ coincide with the cluster algebra and upper cluster algebra of ${\mathcal Q}^-_L$ in the sense of (\cite{bz}). 
\end{Cor}

\medskip

\begin{Cor} For the case of the quantum $n \times r$ matrix algebra,
  the  quantum cluster algebra is equal to its upper bound. This holds irrespective of which $B$ we use in our initial seed.
 \end{Cor}

\end{document}